\documentclass[10pt]{article}

\usepackage{latexsym,amssymb,amsmath,amsthm,paralist}
\def\mapr#1{\stackrel{#1}{\longrightarrow}}

\def\mapd#1{\Big\downarrow\rlap {$\vcenter{\hbox{$\scriptstyle{{#1}}$}}$}}
\def\injr#1{\stackrel{#1}{\lhook\joinrel\relbar\joinrel\rightarrow}}
\def\surjd#1{\lower4pt\hbox{$\downarrow$}\kern-5.65pt\Big\downarrow\rlap {$\vcenter{\hbox{$\scriptstyle{{#1}}$}}$}}
 \def\injd#1{
 \setlength{\unitlength}{0.1pt}
 \begin{picture}(40,120)
 \put(20,102){\vector(0,-1){158}} \put(12.5,100){\mbox{$\scriptscriptstyle
 \cap $}} \put(47,10){\mbox{$\scriptstyle #1 $}}
 \end{picture}}
 \def\eqd{\Big\|}

\newcommand{\Spec}{\text{\it Spec}}
\newcommand{\Q}{{\mathbb Q}}
\newcommand{\Z}{{\mathbb Z}}
\newcommand{\F}{{\mathbb F}}
\newcommand{\Frob}{{\mathit{Frob}}}
\newcommand{\N}{{\mathbb N}}
\newcommand{\G}{{\mathbb G}}
\renewcommand{\O}{{\mathcal{O}}}
\newcommand{\ab}{\text{\it ab}}
\newcommand{\Cl}{\text{\it Cl}}
\newcommand{\p}{{\mathfrak p}}
\newcommand{\q}{{\mathfrak q}}
\newcommand{\sm}{{\,\smallsetminus\,}}
\newcommand{\et}{\text{\it et}}
\newcommand{\fl}{\text{\it fl}\,}
\newcommand{\cd}{\text{\it cd}\:}
\newcommand{\gr}{\textrm{gr}}
\newcommand{\g}{\mathfrak{g}}
\newcommand{\scd}{\text{\it scd}\:}
\newcommand{\Gal}{\text{\it Gal}}

\newcommand{\Hom}{\text{\rm Hom}}
\font\emas = cmsy10 scaled\magstep2
\newcommand{\freeproductmed}{\mathop{\lower.2mm\hbox{\emas \symbol{3}}}\limits}
\newcommand{\lang}{\longrightarrow}
\font\russ=wncyr10
\renewcommand{\min}{\text{\rm min}}
 \newcommand{\ressum}{\mathop{\hbox{${\displaystyle\bigoplus}'$}}\limits}
 \newcommand{\ressumsmall}{\mathop{\hbox{${\bigoplus}'$}}}
 \newcommand{\coker}{\text{\rm coker}}
\newcommand{\ds}{\displaystyle}
\renewcommand{\a}{\mathfrak{a}}
\renewcommand{\r}{\mathfrak{r}}

\newcommand{\nr}{\mathit{nr}}
\renewcommand{\P}{\mathfrak{P}}

\font\russ=wncyr10
 1

\def\Sha{\hbox{\russ\char88}}

\def\Be{\hbox{\russ\char66}}

\newtheoremstyle{alex}% name
  {}%      Space above
  {}%      Space below 20pt plus 5pt minus 3pt
  {\sl}%         Body font
  {}%         Indent amount (empty = no indent, \parindent = para indent)
  {\bf}% Thm head font
  {.}%        Punctuation after thm head
  {.5em}%     Space after thm head: " " = normal interword space;
  {}%         Thm head spec (can be left empty, meaning `normal')
\newtheoremstyle{alexdef}% name
  {}%      Space above
  {}%      Space below
  {\sl}%         Body font
  {}%         Indent amount (empty = no indent, \parindent = para indent)
  {\bf}% Thm head font
  {.}%        Punctuation after thm head
  {.5em}%     Space after thm head: " " = normal interword space;
  {}%         Thm head spec (can be left empty, meaning `normal')
\theoremstyle{alex}
 \newtheorem{itheorem}{Theorem}
\newtheorem{theorem}{Theorem}[section]
\newtheorem{corollary}[theorem]{Corollary}
\newtheorem{lemma}[theorem]{Lemma}
\newtheorem{proposition}[theorem]{Proposition}
\newtheorem{conjecture}{Conjecture}

\theoremstyle{alexdef}
\newtheorem*{definition}{Definition}

\title{\bf\boldmath Rings of integers of type $K(\pi,1)$}
\author{by Alexander Schmidt}
\date{October 17, 2007}

\begin{document}
\maketitle

\begin{abstract} We investigate the Galois group $G_S(p)$ of the maximal $p$-extension unramified outside a finite $S$ of primes of a number field in the (tame) case, when no prime dividing~$p$ is in $S$. We show that the cohomology of $G_S(p)$ is `often' isomorphic to the \'{e}tale cohomology of the scheme $\Spec(\O_k \sm S)$, in particular, $G_S(p)$ is of cohomological dimension~$2$ then.
\end{abstract}

%\begin{minipage}{11cm}
%{\footnotesize
%{\sc\noindent Abstract}:
%We investigate the Galois group $G_S(p)$ of the maximal $p$-extension unramified outside a finite $S$ of primes of a number field in the (tame) case, when no prime dividing $p$ is in $S$. We show that the cohomology of $G_S(p)$ is `often' isomorphic to the \'{e}tale cohomology of the scheme $\Spec(\O_k \sm S)$, in particular, $G_S(p)$ is of cohomological dimension~$2$ then.  }
%\end{minipage}

\section{Introduction}
We call a connected locally noetherian scheme $Y$ a \lq $K(\pi,1)$\rq\ for a prime number~$p$ if the higher homotopy groups of the $p$-completion $Y_\et^{(p)}$ of its \'{e}tale homotopy type $Y_\et$ vanish.  In this paper we consider the case of an arithmetic curve, where the $K(\pi,1)$-property is linked with open questions in the theory of Galois groups with restricted ramification of number fields:

\smallskip
Let $k$ be a number field, $S$ a finite set of nonarchimedean primes of $k$ and~$p$\/ a prime number.  For simplicity, we assume that $p$ is odd or that $k$ is totally imaginary. By a $p$-extension we understand a Galois extension whose Galois group is a \text{(pro-)} $p$-group. Let $k_S(p)$ denote the maximal $p$-extension of $k$ unramified outside $S$ and  put $G_S(p)=\Gal(k_S(p)|k)$.
A systematic study of this group had been started by \v Safarevi\v c, and was continued by  Koch, Kuz'min, Wingberg and many others;  see \cite{NSW}, VIII,\,\S 7 for basic properties of $G_S(p)$. In geometric terms (and omitting the base point), we have
\[
G_S(p)\cong \pi_1\big((\Spec(\O_k)\sm S)_\et^{(p)}\big).
\]
As is well known to the experts, if $S$ contains the set $S_p$ of primes dividing~$p$, then $\Spec(\O_k)\sm S$ is a $K(\pi,1)$ for~$p$ (see Proposition~\ref{wildiskpi1} below). In particular, if $S\supset S_p$, then $G_S(p)$ is of cohomological dimension less or equal to $2$.

\medskip
The group $G_S(p)$  is most mysterious in the {\it tame} case, i.e.\ when $S\cap S_p=\varnothing$. In this case, examples when $\Spec(\O_k)\sm S$ is {\it not} a $K(\pi,1)$ are easily constructed. On the contrary, until recently not a single $K(\pi,1)$-example was known.   The following properties of the group $G_S(p)$ were known so far:

\pagebreak
%\smallskip
\begin{compactitem}
\item $G_S(p)$ is a `fab-group', i.e.\ $U^\ab$ is finite for each open subgroup $U\subset G$.
\item $G_S(p)$ can be infinite (Golod-\v Safarevi\v c).
\item $G_S(p)$ is a finitely presented pro-$p$-group (Koch).
\end{compactitem}

\smallskip\noindent
A conjecture of Fontaine and Mazur \cite{FM} asserts that $G_S(p)$ has no infinite $p$-adic analytic quotients.

\bigskip
In 2005, Labute considered the case $k=\Q$ and found finite sets  $S$ of prime numbers (called strictly circular sets) with $p\notin S$ such that $G_S(p)$ has cohomological dimension~$2$. In \cite{S} the author showed that, in the examples given by Labute,  $\Spec(\Z)\sm S$ is a $K(\pi,1)$ for $p$.

The objective of this paper is a systematic study of the $K(\pi,1)$-property. Our focus is on the tame case, where we conjecture that rings of integers of type $K(\pi,1)$ are cofinal in the following sense:

\begin{conjecture}\label{conj}
Let $k$ be a number field and let $p$ be a prime number. Assume that $p\neq 2$ or that $k$ is totally imaginary. Let $S$ be a finite set of primes of $k$ with $S\cap S_p=\varnothing$. Let, in addition, a set $T$ of primes of Dirichlet density $\delta(T)=1$ be given. Then there exists a finite subset $T_1\subset T$ such that
$
\Spec(\O_k)\sm (S\cup T_1)$
is a $K(\pi,1)$ for $p$.
\end{conjecture}
Of course we may assume that $T\cap S_p=\varnothing$ in the conjecture. Our main result is the following

\begin{itheorem} \label{exthmi}
Conjecture~\ref{conj} is true if the number field $k$ does not contain a primitive $p$-th root of unity and the class number of\/ $k$ is prime to $p$.
\end{itheorem}

Explicit examples of rings of integers of type $K(\pi,1)$ can be found in \cite{La}, \cite{S} (for $k=\Q$) and in \cite{Vo} (for $k$~imaginary quadratic).

\medskip
The $K(\pi,1)$-property has strong consequences. We write $X=\Spec(\O_k)$ and assume in all results that $p\neq 2$ or that $k$ is totally imaginary.
Primes $\p\in S\sm S_p$ with $\mu_p\not\subset k_\p$ are redundant in the sense that removing these primes from~$S$ does not change $(X\sm S)_\et^{(p)}$, see section~\ref{sect3}. In the tame case, we may therefore restrict our considerations to sets of primes whose norms are congruent to $1$ modulo $p$. These are the results.

\begin{itheorem} \label{1.1} Let $S$ be a finite non-empty set of primes of\/ $k$ whose norms are congruent to $1$ modulo $p$.  If\/ $X\sm S$ is a $K(\pi,1)$ for $p$ and $G_S(p)\neq 1$, then the following hold.
\begin{itemize}
\item[\rm (i)] $\cd G_S(p)=2$, $\scd G_S(p)=3$.
\item[\rm (ii)] $G_S(p)$ is a duality group.
\end{itemize}
The dualizing module $D$ of $G_S(p)$ is given by
$
D=\text{tor}_p C_S(k_S(p))$, i.e.\ it is the subgroup of $p$-torsion elements in the $S$-id\`{e}le class group of\/ $k_S(p)$.
\end{itheorem}

\noindent
{\bf Remarks:} 1. If $X\sm S$ is a $K(\pi,1)$ for $p$ and $G_S(p)= 1$, then  $k$ is imaginary qua\-dratic, $\#S=1$ and $p=2$ or $3$. See Proposition~\ref{degen} for a more precise statement.

\smallskip\noindent
2. We have a natural exact sequence

\bigskip
$\displaystyle 0 \to \mu_{p^\infty}(k_S(p)) \to \bigoplus_{{\mathfrak p} \in S} \text{\rm Ind\,}^{G_{\mathfrak p}(k_S(p)|k)}_{\Gal(k_S(p)|k)}\ \mu_{p^\infty}(k_S(p)_\p) \to
$
\begin{flushright}
$\text{tor}_p C_S(k_S(p)) \to {\mathcal O}_{k_S(p), S}^\times \otimes \Q_p/\Z_p \to 0\ ,\quad $
\end{flushright}
where ${\mathcal O}_{k_S(p), S}^\times$ is the group of $S$-units of\/ $k_S(p)$ and $\mu_{p^\infty}(K)$ denotes the group of all $p$-power roots of unity in a field $K$. Note that $\mu_{p^\infty}(k_S(p))$ is finite, while, by Theorem~\ref{1.2} below,  for $\p\in S$ the field $k_S(p)_\p$ contains all $p$-power roots of unity.

\smallskip\noindent
3. In the wild case $S\supset S_p$, where $X\sm S$ is always a $K(\pi,1)$ for $p$, $G_S(p)$ is of cohomological dimension~$1$ or $2$.  The strict cohomological dimension is conjecturally equal to $2$ (=Leopoldt's conjecture for each finite subextension of $k$ in $k_S(p)$). In the wild case, $G_S(p)$ is often, but not always a duality group, cf.\ \cite{NSW} Prop.\ 10.7.13.

\bigskip
Allowing ramification at a prime $\p$ does not mean that the ramification is realized globally.
Therefore it is a natural and interesting question how far we get locally at the primes in $S$ when going up to $k_S(p)$.  See \cite{NSW} X,\,\S 3 for results in the wild case. In the tame case, we have the following
\begin{itheorem} \label{1.2} Let $S$ be a finite non-empty set of primes of\/ $k$ whose norms are congruent to\/ $1$ modulo $p$.  If\/ $X\sm S$ is a $K(\pi,1)$ for\/  $p$ and $G_S(p)\neq 1$, then
\[
k_S(p)_\p=k_\p(p)
\]
for all primes $\p\in S$, i.e.\ $k_S(p)$ realizes the maximal $p$-extension of the local field~$k_\p$.
\end{itheorem}

\noindent
{\bf Remark:} Under the assumptions of the theorem, let $\q\notin S$. Then either $\q$ splits completely in $k_S(p)$, or $k_S(p)$ realizes the maximal unramified $p$-extension $k_\q^\nr(p)$. We do not know whether the completely split case actually occurs.

\bigskip
The next result addresses the question of enlarging the set $S$ without destroying the $K(\pi,1)$-property.

\begin{itheorem} \label{1.3} Let $S'$ be a finite non-empty set of primes of $k$ whose norms are congruent to $1$ modulo $p$ and let $S\subset S'$ be a nonempty subset. Assume that $X\sm S$ is a $K(\pi,1)$ for $p$ and that $G_S(p)\neq 1$. If each $\q\in S'\sm S$ does not split completely in $k_S(p)$, then $X\sm S'$ is a $K(\pi,1)$ for $p$. Furthermore, in this case, the arithmetic form of Riemann's existence theorem holds:  the natural homomorphism
\[
\freeproductmed_{{\mathfrak p} \in S'\backslash S(k_S(p))} T_\p(k_{S'}(p)|k_S(p)) \longrightarrow \Gal(k_{S'}(p)|k_S(p))
\]
is an isomorphism, i.e.\ $\Gal(k_{S'}(p)|k_S(p))$ is the free pro-$p$ product of a bundle of inertia groups.
\end{itheorem}

Finally, we deduce a statement on universal norms of unit groups.

\begin{itheorem} \label{1.4} Let $S$ be a finite non-empty set of primes of\/ $k$ whose norms are congruent to $1$ modulo $p$. Assume that $X\sm S$ is a $K(\pi,1)$ for $p$ and that $G_S(p)\neq 1$. Then
\[
\varprojlim_{K\subset k_S(p)} \O_K^\times \otimes \Z_p\,=\,0\,=\,\varprojlim_{K\subset k_S(p)} \O_{K,S}^\times \otimes \Z_p,
\]
where $K$ runs through all finite subextensions of\/ $k$ in $k_S(p)$, $\O_K^\times$ and $\O_{K,S}^\times$ are the groups of units and $S$-units of\/ $K$, respectively, and the transition maps are the norm maps.
\end{itheorem}

\bigskip
The structure of this paper is as follows. First we give the necessary definitions and make some calculations of \'{e}tale cohomology groups for which we couldn't find an appropriate reference. In section~\ref{sect3}, we deal with the first obstruction against the $K(\pi,1)$-property, the $h^2$-defect. Then we recall Labute's results on mild pro-$p$-groups, which we use in the proof of Theorem~\ref{exthmi} given in section~6. In the last three sections we prove Theorems 2--5.

\bigskip
The author thanks Denis Vogel for valuable comments on a preliminary version of this paper.

\section{First observations}\label{1}

We tacitly assume all schemes to be connected and omit base points from the notation. Let $Y$ be a locally noetherian scheme and let $p$ be a prime number. We denote by $Y_\et^{(p)}$ the $p$-completion of the \'{e}tale homotopy type of $Y$, see \cite{AM}, \cite{Fr}. By $\widetilde{Y}{(p)}$ we denote the universal pro-$p$-covering of $Y$. The projection $\widetilde{Y}{(p)}\to Y$ is Galois with group $\pi_1^\et(Y)(p)=\pi_1(Y_\et^{(p)})$, cf.\ \cite{AM}, (3.7). Any discrete $p$-torsion $\pi_1^\et(Y)(p)$-module $M$ defines a locally constant sheaf on $Y_\et$, which we denote by the same letter.
The Hochschild-Serre spectral sequence defines natural homomorphisms
\[
\phi_{M,i}: H^i(\pi_1^\et(Y)(p), M) \longrightarrow H^i_\et(Y, M), \ i\geq 0.
\]
Since  $H^1_\et(\widetilde{Y}{(p)}, M)=0$, the map $\phi_{M,i}$ is an isomorphism for $i=0$ and $1$, and is injective for $i=2$.
For a pro-$p$-group $G$ we denote by $K(G,1)$ the associated Eilenberg-MacLane space (\cite{AM}, (2.6)).
\begin{proposition} \label{kpi1cond}
The following conditions are equivalent:
\begin{itemize}
\item[\rm (i)] The classifying map $Y_\et^{(p)} \longrightarrow K(\pi_1^\et(Y)(p),1)$ is a weak equivalence.
\item[\rm (ii)] $\pi_i(Y_\et^{(p)})=0$ for all $i\geq 2$.
\item[\rm (iii)] $H^i_\et (\widetilde{Y}{(p)}, \Z/p\Z)=0$ for all $i\geq 1$.
\item[\rm (iv)] $\phi_{\Z/p\Z,i}: H^i(\pi_1^\et(Y)(p),  \Z/p\Z) \longrightarrow H^i_\et(Y,  \Z/p\Z)$ is an isomorphism for all $i\geq 0$.
\item[\rm (v)] $\phi_{M,i}: H^i(\pi_1^\et(Y)(p), M) \longrightarrow H^i_\et(Y, M)$ is an isomorphism for all $i\geq 0$ and  any discrete $p$-torsion $\pi_1^\et(Y)(p)$-module $M$.
\end{itemize}
\end{proposition}

\begin{proof}
The equivalences (i)$\Leftrightarrow$(ii)$\Leftrightarrow$(v) are the content of  \cite{AM}, (4.3), (4.4). The equivalence (iii)$\Leftrightarrow$(iv) follows in a straightforward manner from the Hochschild-Serre spectral sequence. The implication (v)$\Rightarrow$(iv) is trivial.

Assume that (iv) holds. As $\pi_1^\et(Y)(p)$ is a pro-$p$-group, any finite $p$-primary $\pi_1^\et(Y)(p)$-module $M$ has a composition series with graded pieces isomorphic to $\Z/p\Z$ with trivial $\pi_1^\et(Y)(p)$-action (\cite{NSW}, Corollary 1.7.4). Hence, if $M$ is finite, the five-lemma implies that $\phi_{M, i}$ is an isomorphism for all $i$.  An arbitrary discrete $p$-primary $\pi_1^\et(Y)(p)$-module is the filtered inductive limit of finite $p$-primary $\pi_1^\et(Y)(p)$-modules.  Since group cohomology (\cite{NSW}, Proposition 1.5.1) and \'{e}tale cohomology (\cite{AGV}, VII, 3.3) commute with filtered inductive limits, $\phi_{M, i}$ is an isomorphism  for all $i$ and all discrete $p$-torsion $\pi_1^\et(Y)(p)$-modules $M$. This implies (v) and completes the proof.
\end{proof}

\begin{definition}
We say that $Y$ is a {\bf\boldmath $K(\pi,1)$ for\/ $p$} if the equivalent conditions of Proposition~\ref{kpi1cond} are satisfied.
\end{definition}

Now let $k$ be a number field, $S$ a finite set of nonarchimedean primes of $k$ and~$p$\/ a prime number. We put $X=\Spec(\O_k)$.  The following observation is straightforward.

\begin{corollary}\label{finite_ext}
Let $k'$ be a finite subextension of\/ $k$ in $k_S(p)$ and let $X'=\Spec(\O_{k'})$, $S'=S(k')$. Then the following are equivalent.
\begin{itemize}
\item[\rm (i)] $X\sm S$ is a $K(\pi,1)$ for $p$,
\item[\rm (ii)] $X'\sm S'$ is a $K(\pi,1)$ for $p$.
\end{itemize}
\end{corollary}

\begin{proof}
Both schemes have the same universal pro-$p$-covering.
\end{proof}

\bigskip
We denote by $S_p$ and $S_\infty$ the set of primes of $k$ dividing $p$ and the set of archimedean primes of $k$, respectively. For a set $S$ of primes (which may contain archimedean places), let $k_S(p)$ be the maximal $p$-extension of $k$ unramified outside $S$ and $G_S(p)=\Gal(k_S(p)|k)$. For a finite set $S$ of nonarchimedean primes of $k$ we have the identification
\[
\pi_1^\et((X\sm S)_\et^{(p)})= G_{S\cup S_\infty}(p).
\]
If $p$ is odd or $k$ is totally imaginary, then $G_S(p)=G_{S\cup S_\infty}(p)$.
The following proposition is given for sake of completeness. It deals with the `wild' case $S\supset S_p$, and is well known.
\begin{proposition} \label{wildiskpi1}
If\/ $S$ contains $S_p$, then $X\sm S$ is a $K(\pi,1)$ for $p$.
\end{proposition}
\begin{proof}  We verify condition (v) of Proposition~\ref{kpi1cond}. Let $k_{S\cup S_\infty}$ be the maximal extension of $k$ unramified outside $S\cup S_\infty$ and put $G_{S\cup S_\infty}=Gal(k_{S\cup S_\infty}|k)$. For any $p$-primary discrete $G_{S\cup S_\infty}(p)$-module $M$ the homomorphism $\phi_{M,i}$ factors as
\[
H^i(G_{S\cup S_\infty}(p),M) \to H^i(G_{S\cup S_\infty},M) \to H^i_\et(X\sm S, M).
\]
By \cite{NSW}, Cor.~10.4.8, the left map is an isomorphism. That also the right map is an isomorphism follows in a straightforward manner by using the Kummer sequence, the Principal Ideal Theorem and known properties of the Brauer group, see for example \cite{Zi}, Prop. 3.3.1.\ or \cite{Mi}, II Prop.~2.9.
\end{proof}

\noindent
{\bf Remark:} If $p=2$ and $k$ has real places it is useful to work with the modified \'{e}tale site defined by T.~Zink \cite{Zi}, which takes the real archimedean places of~$k$ into account. Proposition~\ref{wildiskpi1} holds true also for the modified \'{e}tale site,  see \cite{twoinf}, Thm.~6.

\section{Calculation of \'{e}tale cohomology groups}\label{calcsect}

As a basis of our investigations, we need the calculation of the \'{e}tale cohomology groups of open subschemes of $\Spec(\O_k)$ with values in the constant sheaf $\Z/p\Z$.
Let $p$ be a fixed prime number. All cohomology groups are taken with respect to the constant sheaf $\Z/p\Z$, which we omit from the notation. Furthermore, we use the notation
\[
h^i(-)=\dim_{\F_p} H^i_\et(-)\quad(=\dim_{\F_p} H^i_\et(-,\Z/p\Z)\;)
\]
for the occurring cohomology groups.  We start with some well-known local computations.

\begin{proposition}\label{localcoh}
Let $k$ be a nonarchimedean local field of characteristic zero and residue characteristic $\ell$. Let $X=\Spec(\O_k)$ and  let $x$ be the closed point of~$X$.   Then the \'{e}tale local cohomology groups $H^i_x(X)$ vanish for $i\leq 1$ and $i\geq 4$, and
\[
h^2_x(X)= \left\{
\begin{array}{cl}
\delta,& \text{ if }\ \ell\neq p,\\
\delta+[k:\Q_p],& \text{ if }\ \ell=p,
\end{array}\right.
\]
where $\delta=1$ if\/ $\mu_p\subset k$ and zero otherwise.
Furthermore, $h^3_x(X)=\delta$. In particular, we have the Euler-Poincar\'{e} characteristic formula
\[
\sum_{i=0}^3 (-1)^i h^i_x(X)= \left\{
\begin{array}{cl}
0,& \text{ if }\ \ell\neq p,\\
\,[k:\Q_p],& \text{ if }\ \ell=p.
\end{array}\right.
\]
\end{proposition}

\begin{proof}
As $X$ is henselian, we have isomorphisms $H^i_\et(X) \cong H^i_\et(x)$ for all $i$, and therefore
\[
h^i(X)= \left\{
\begin{array}{cl}
1 & \hbox{ for } i=0,1,\\
0 & \hbox{ for } i\geq 2.
\end{array}
\right.
\]
Furthermore, $X\sm \{x\}=\Spec(k)$, hence $H^i_\et(X\sm \{x\})\cong H^i(k)$. The local duality theorem (cf.\ \cite{NSW}, Theorem 7.2.15) shows $h^2(X\sm \{ x\})=\delta$, and by \cite{NSW}, Corollary 7.3.9, we have
\[
h^1(X\sm \{ x\})= \left\{
\begin{array}{cl}
1+\delta & \hbox{ if } \ell\neq p,\\
1+\delta+[k:\Q_p] & \hbox{ if } \ell=p.
\end{array}
\right.
\]
Furthermore, the natural homomorphism $H^1_\et(X) \to H^1_\et(X\sm \{ x\})$ is injective. Therefore the result of the proposition follows from the exact excision sequence
\[
\cdots \to H^i_x(X) \to H^i_\et(X) \to H^i_\et(X\sm \{x\}) \to H^{i+1}_x(X) \to \cdots\;.
\]
\end{proof}

Now let $k$ be a number field, $S$ a finite set of nonarchimedean primes of $k$ and $X=\Spec(\O_k)$.  We assume for simplicity that $p$ is odd or that $k$ is totally imaginary, so that we can ignore the archimedean places of $k$ for cohomological considerations. We introduce the following notation
\begin{tabbing}
\quad \= $r_1$\hspace*{1cm} \= the number of real places of $k$\\
\>$r_2$\ \> the number of complex places of $k$\\
\>$r$\ \> $=r_1+r_2$, the number of  archimedean places of $k$\\
\> $S_p$ \> the set of places of $k$ dividing $p$\\
\> $\delta$ \> equal to $1$ if $\mu_p\subset k$ and zero otherwise\\
\> $\delta_\p$ \> equal to $1$ if $\mu_p\subset k_\p$ and zero otherwise\\
\> $\Cl(k)$\> the ideal class group of $k$\\
\> $\Cl_S(k)$\> the $S$-ideal class group of $k$\\
\> $h_k$\> $=\# \Cl(k)$, the class number of $k$\\
\> $_n A$\> $=\ker(A \stackrel{\cdot n}{\to} A)$, where $A$ is an abelian group and $n\in \N$ \\
\> $A/n$\> $=\coker(A \stackrel{\cdot n}{\to} A)$, where $A$ is an abelian group and $n\in \N$.
\end{tabbing}

\begin{proposition}\label{globalchi}
Assume that $p\neq 2$ or that $k$ is totally imaginary.
Then $H^i_\et(X\sm S)=0$ for\/ $i\geq 4$, and
\[
\chi(X\sm S):= \sum_{i=0}^3 (-1)^i h^i(X\sm S)= r - \sum_{\p \in S\cap S_p} [k_\p:\Q_p]\;.
\]
In particular,
\[
\chi(X\sm S)=\left\{
\begin{array}{cl}
r, & \hbox{ if }\ S\cap S_p=\varnothing\,,\\
-r_2, & \hbox{ if }\ S\supset S_p\, .
\end{array}
\right.
\]
\end{proposition}

\begin{proof}
The assertion for $S=S_p$ is well known, see \cite{Mi}, II Theorem~2.13 (a). Consider the exact excision sequence
\[
\cdots \to \bigoplus_{\p \in S} H^i_\p(X_\p) \to H^i_\et(X) \to H^i_\et (X\sm S) \to \bigoplus_{\p \in S} H^{i+1}_\p(X_\p) \to \cdots \ ,
\]
where $X_\p=\Spec(\O_{k,\p})$ is the spectrum of the completion of $\O_k$ at $\p$. Using this excision sequence for $S=S_p$, Proposition~\ref{localcoh} implies the result for $S=\varnothing$, noting that $\sum_{\p \in S_p}[k_\p:\Q_p]-r_2=[k:\Q]-r_2=r$. The result for arbitrary $S$ follows from the case $S=\varnothing$, the above excision sequence and from Proposition~\ref{localcoh}.
\end{proof}

In order to give formulae for the individual cohomology groups, we consider the Kummer group
 (cf.\ \cite{NSW}, VIII,\, \S 6)
\[
V_S(k):=\{ a \in k^\times \mid a \in k_\p^{\times p} \hbox{ for } \p \in S \hbox{ and } a \in U_\p k_\p^{\times p} \hbox{ for } \p \notin S\}/k^{\times p},
\]
where $U_\p$ denotes the unit group of the local field $k_\p$ (convention: $U_\p=k_\p^\times$ if $\p$ is archimedean).\footnote{In terms of flat cohomology, we have
$
V_S(k)= \ker \big(H^1_\fl(X\sm S,\mu_p) \to \bigoplus_{\p \in S} H^1(k_\p,\mu_p)\big).
$}
$V_S(k)$ is a finite group. More precisely, we have the following

\begin{proposition}\label{VSchange}
There exists a natural exact sequence
\[
0\longrightarrow \O_k^\times /p \longrightarrow V_\varnothing(k) \longrightarrow \null_p \Cl(k) \longrightarrow 0\,.
\]
In particular,
\[
\dim_{\F_p} V_\varnothing(k) = \dim_{\F_p} \null_p \Cl(k) + \dim_{\F_p} \O_k^\times /p= \dim_{\F_p} \null_p \Cl(k) + r-1+\delta.
\]
If\/ $S$ is arbitrary and \/ $\p\notin S$ is an additional prime of\/ $k$, then we have a natural exact sequence
\[
0 \longrightarrow V_{S \cup \{\p\}}(k)\stackrel{\phi}{\longrightarrow} V_S(k)  \longrightarrow U_\p k_\p^{\times p}/k_\p^{\times p}.
\]
For $\p\notin S_p$, we have
$
\dim_{\F_p} \coker(\phi)\leq \delta_\p$.
\end{proposition}

\begin{proof}
Sending an $a\in V_\varnothing(k)$ to the class in $\Cl(k)$ of the fractional ideal $\mathfrak a$ with $(a)={\mathfrak a}^p$ yields a surjective homomorphism $V_\varnothing(k) \to \null_p\Cl(k)$ with kernel
$\O_k^\times /p$. This, together with Dirichlet's Unit Theorem, shows the first statement.
The second exact sequence follows immediately from the definitions.
There are natural isomorphisms
\[
U_\p k_\p^{\times p}/k_\p^{\times p}\cong U_\p/U_\p\cap k_\p^{\times p} = U_\p/ U_\p^{ p}.
\]
For $\p\notin S_p$ we have $\dim_{\F_p} U_\p/ U_\p^{p} = \delta_\p$,  showing the last statement.
\end{proof}

The \'{e}tale cohomology groups of $X\sm S$ have the following dimensions.
\renewcommand{\arraystretch}{1.5}
\begin{theorem} \label{globcoh} Let $S$ be a finite set of nonarchimedean primes of\/ $k$. Assume $p\neq 2$ or that $k$ is totally imaginary. Then $H^i_\et(X\sm S)=0$ for $i\geq 4$ and
\[
\begin{array}{lcl}
h^0(X \sm S)&=& 1\, ,\\
h^1(X \sm S) & = & 1+  \ds\sum_{\p\in S} \delta_\p
               - \delta + \dim_{\F_p}V_S +  \ds\sum_{\p\in S\cap S_p} [k_\p:\Q_p] -r,\\
h^2(X \sm S)& =& \ds\sum_{\p\in S} \delta_\p
               - \delta + \dim_{\F_p}V_S +\theta \, ,\\
h^3(X \sm S)&=& \theta\, .\\
\end{array}
\]
Here $\theta$ is equal to $1$ if\/ $\delta =1$ and $S=\varnothing$, and zero in all other cases.
\end{theorem}
\renewcommand{\arraystretch}{1}

\begin{proof}
The statement on $h^0$ is trivial and the vanishing of $H^i$ for $i\geq 4$ was already part of Proposition~\ref{globalchi}.  Artin-Verdier duality (see \cite{Ma}, 2.4 or \cite{Mi}, Theorem~3.1) shows
\[
H^3_\et(X)^\vee\cong\Hom_X(\Z/p\Z, \G_m)=\mu_p(k).
\]
Consider the exact excision sequence
\[
\bigoplus_{\p\in S} H^3_\p(X_\p) \stackrel{\alpha}{\to} H^3_\et (X) \stackrel{\beta}{\to} H^3_\et(X\sm S) {\to} \bigoplus_{\p\in S} H^4_\p(X_\p).
\]
By Proposition~\ref{localcoh}, the right hand group is zero, hence $\beta$ is surjective.  By the local duality theorem (see \cite{Ma}, 2.4, \cite{Mi}, II Corollary 1.10), the dual map  to $\alpha$ is the natural inclusion
\[
 \mu_p(k) \to \bigoplus_{\p\in S} \mu_p(k_\p),
\]
which is injective, unless $\delta=1$ and $S=\varnothing$.
Therefore $h^3(X \sm S)=1$ if $\delta=1$ and $S = \varnothing$, and zero otherwise.
Using the isomorphism
$H^1(G_S(p)) \stackrel{\sim}{\rightarrow} H^1_\et(X\sm S)$, the statement on $h^1$ follows from the corresponding formula for the first cohomology of $G_S(p)$, see \cite{NSW}, Theorem 8.7.11. Finally, the result for $h^2$ follows by using the Euler-Poincar\'{e} characteristic formula in Proposition~\ref{globalchi}.
\end{proof}

\begin{corollary} \label{condforkpi1}
Assume that $\delta=0$ or $S\neq\varnothing$. Then  $X\sm S$ is a $K(\pi,1)$ for $p$ if and only if the following conditions {\rm (i)} and {\rm (ii)} are satisfied.
\begin{itemize}
\item[\rm (i)] $\phi_{2}\colon\ H^2(G_S(p))\hookrightarrow H^2_\et(X\sm S)$ is an isomorphism,
\item[\rm (ii)] $\cd G_S(p) \leq 2$.
\end{itemize}
\end{corollary}

\begin{proof} The given conditions are obviously necessary.
Furthermore,  $\phi_{0}$  and $\phi_{1}$ are isomorphisms and $H^i_\et(X\sm S)=0$ for $i\geq 3$ by Theorem~\ref{globcoh}. Therefore (i) and (ii) imply that $\phi_{i}$ is an isomorphism for all $i$. Hence condition (iv) of Proposition~\ref{kpi1cond} is satisfied for $X\sm S$ and $p$.
\end{proof}

\bigskip
Let $F$ be a  locally constant sheaf on $(X\sm S)_\et$.
For each prime $\p$ the composite map $\O_{k,S} \to k \to k_\p$ induces natural maps $H^i_\et (X\sm S, F)\to H^i(k_\p, F)$ for all $i\geq 0$.

\begin{definition} For any locally constant sheaf $F$ on $(X\sm S)_\et$ we put
\[
\Sha^i(k,S,F):= \ker \Big( H^i_\et (X\sm S, F)\longrightarrow \prod_{\p\in S} H^i(k_\p, F) \Big).
\]
Assume a prime number $p$ is fixed. Then we write $\Sha^i(k,S):=\Sha^i(k,S,\Z/p\Z)$ and,  following historical notation, we put $\Be_S(k):= V_S(k)^\vee$, where
 $\scriptstyle \vee$ denotes the Pontryagin dual.
\end{definition}

The next theorem is sharper than \cite{NSW}, Thm.~8.7.4, as the group $\Sha^2(G_S)$, which was considered there, is a subgroup of $\Sha^2(k,S)$.   If $p=2$ and $k$ has real places, then Theorem~\ref{shaincl} remains true after replacing \'{e}tale cohomology with its modified version.

\begin{theorem} \label{shaincl} Assume $p\neq 2$ or that $k$ is totally imaginary. Then
there exists a natural isomorphism
\[
\Sha^2(k,S) \mapr{\sim} \Be_S(k).
\]
\end{theorem}

\begin{proof} The proof of \cite{NSW}, Thm.~8.7.4  can be adapted to show also the stronger statement here. However, we decided to take the short cut by using flat duality.
For any prime $\p$ of $k$ one easily computes the local cohomology groups for the flat topology as  $H^1_{\fl, \p}(X,\mu_p)=0$  and $H^2_{\fl, \p}(X,\mu_p) \cong k_\p^\times/U_\p k_\p^{\times p}$.
Therefore excision and Kummer theory imply an exact sequence
\[
0\to H^1_\fl(X,\mu_p) \to k^\times/k^{\times p} \to \bigoplus_\p
k_\p^\times/U_\p k_\p^{\times p}.
\]
As $H^1_\fl(X_\p^h,\mu_p)\cong U_\p/p$, we obtain the exact sequence
\[
0 \to V_S(k) \to H^1_\fl(X,\mu_p) \to \bigoplus_{\p \in S} H^1_{\fl} (X_\p^h,\mu_p).\leqno (*)
\]
By excision, and noting $H^3_\p(X,\Z/p\Z)\cong H^2(k_\p, \Z/p\Z)$, we have an exact sequence
\[
\bigoplus_{\p\in S} H^2_\p(X, \Z/p\Z) \to H^2_\et(X, \Z/p\Z) \to \Sha^2(k,S) \to 0. \leqno(**)
\]
Comparing sequences $(*)$ and $(**)$ via local and global flat duality, we obtain the asserted isomorphism.
\end{proof}

We provide the following lemma for further use.
\begin{lemma}\label{h3oben}
Let $K\subset k_S(p)$ be an extension of\/ $k$ inside $k_S(p)$ and let $(X\sm S)_K$ be the normalization of\/ $X\sm S$ in $K$. If\/ $\delta=0$, or $S\neq \varnothing$ or  $K|k$ is infinite,  then
\[
H^3_\et ((X\sm S)_K)=0.
\]
\end{lemma}

\begin{proof}
We denote the normalization of $X\sm S$ in any algebraic extension field $k'$ of $k$ by $(X\sm S)_{k'}$.
\'{E}tale cohomology commutes with inverse limits of schemes if the transition maps are affine (see \cite{AGV}, VII, 5.8). Therefore
\[
H^3((X\sm S)_K)=\varinjlim_{k'\subset K} H^3((X\sm S)_{k'}),
\]
where $k'$ runs through all finite subextensions of\/ $k$ in $K$.
If $\delta=0$ or $S\neq \varnothing$, then, by Theorem~\ref{globcoh}, $H^3_\et((X\sm S)_{k'})=0$ for all $k'$ and the limit is obviously zero. Assume $\delta=1$ and $S=\varnothing$.
Then, by Artin-Verdier duality,
\[
H^3_\et(X_{k'}) \cong \mu_p(k')^\vee.
\]
For $k' \subset k'' \subset K$, the transition map
\[
H^3_\et(X_{k'}) \to H^3_\et(X_{k''})
\]
is the dual of the norm map $N_{k''|k'}\colon\ \mu_p(k'') \to \mu_p(k')$, hence the zero map if $k' \neq k''$. As $K|k$ is infinite, the limit vanishes.
\end{proof}

\section{\boldmath Removing the $h^2$-defect}\label{sect3}

We start by extending the notions introduced before to infinite sets of primes~$S$. Let $k$ be a number field and $S$ a set of nonarchimedean primes of $k$. We set $X=\Spec(\O_k)$ and
\[
X \sm S =\Spec(\O_{k,S}),
\]
which makes sense also if $S$ is infinite. Let $F$ be a sheaf on $X\sm S$ which comes by restriction from $X\sm T$ for some finite subset $T\subset S$. As each open subscheme of $X$ is affine, we have
\[
H^i_\et(X\sm S, F)  \cong \varinjlim_{\substack{T\subset S'\subset S\\ S'\ \mathrm{ finite}}} H^i_\et(X\sm S', F)
\]
for all $i\geq 0$.

\medskip
We fix a prime number $p$ and put the running assumption that $k$ is totally imaginary if $p=2$. Hence we may ignore archimedean primes for cohomological considerations. The notion of being a $K(\pi,1)$ for $p$ extends in an obvious manner to the case when $S$ is infinite.
Also the isomorphism
\[
\Sha^2(k,S) \mapr{\sim} \Be_S(k)
\]
generalizes to infinite $S$ by passing to the limit over all finite subsets $S'\subset S$.
In particular, $\Sha^2(k,S)$ is finite.

\bigskip
For the remainder of this paper, we assume that $S\cap S_p=\varnothing$. We also keep the running assumption $p\neq 2$ or $k$ is totally imaginary.
\bigskip

\begin{minipage}{11.5cm}
{\it  For shorter notation, we drop $p$ wherever possible. We write  $G_S$ instead of\/ $G_S(p)$,  $k_S$ for $k_S(p)$, and so on.  Unless mentioned otherwise, all cohomology groups are taken with values in $\Z/p\Z$. We keep this notational convention for the rest of this paper.
}
\end{minipage}

\bigskip\noindent
If $\p \nmid p$ is a prime with $\mu_p\not\subset k_\p$, then every $p$-extension of the local field~$k_\p$ is unra\-mi\-fied (see \cite{NSW}, Proposition 7.5.1).  Therefore primes $\p \notin  S_p$ with $N(\p) \not \equiv 1 \bmod p$  cannot ramify in a $p$-extension.
Removing all these redundant primes from $S$, we obtain a subset
$S_{\min} \subset S$ which has the property that $G_S=G_{S_{\min}}$. Moreover, we have the
\begin{lemma}\label{smin}
The natural map
$$
(X\sm S)_\et^{(p)} \longrightarrow (X\sm S_\min)_\et^{(p)}
$$
is a weak homotopy equivalence.
\end{lemma}
\begin{proof}
By \cite{AM}, (4.3), it suffices to show that for every discrete $p$-primary $G_S$-module $M$ the natural maps $H^i_\et(X\sm S_\min,M) \to H^i_\et(X\sm S,M)$ are isomorphisms for all $i$. By the same argument, as in the proof of Proposition~\ref{kpi1cond}, (iv)$\Rightarrow$(v), we may suppose that $M=\Z/p\Z$. Using the excision sequence, it therefore suffices to show that the group $H^i_\p(X\sm S_\min,\Z/p\Z)$ vanishes for all $\p\in S\sm S_\min$. This follows from Proposition~\ref{localcoh}.
\end{proof}

Therefore we can replace $S$ by  $S_\min$ and make the following notational convention for the rest of this paper.

\bigskip\noindent

\begin{minipage}{11cm}{\it  The word `prime' means nonarchimedean prime with norm $\equiv 1 \bmod p$.}
\end{minipage}

\bigskip\noindent
At this point it is useful to redefine the notion of Dirichlet density.

\begin{definition}
Let $S$ be a set of primes of $k$ (of norm $\equiv 1 \bmod p$). The $p$-density $\varDelta^p(S)$ is defined by
\[
\varDelta^p(S) = \delta_{k(\mu_p)}\big(S(k(\mu_p))\big),
\]
where $S(k(\mu_p))$ is the set of prolongations of primes in $S$ to $k(\mu_p)$ and $\delta_{k(\mu_p)}$ denotes the Dirichlet density on the level $k(\mu_p)$. In other words,
\[
\varDelta^p(S)= d\cdot \delta_k(S),\ \text{where}\ \, d=[k(\mu_p):k].
\]
\end{definition}
The set of all primes (of norm $\equiv 1 \bmod p$) has $p$-density equal to~$1$.

\begin{proposition}\label{sha2gleich0}
Let  $S$ be a set of primes of $p$-density $\varDelta^p(S)=1$.  Then there exists a finite subset $T\subset S$ with $\Be_T(k)=0$. In particular, $\Be_S(k)=0= \Sha^2(k,S)$.
\end{proposition}

\begin{proof}
By the Hasse principle for the module $\mu_p$, see \cite{NSW}, Thm.~9.1.3\,(ii), and Kummer theory, the natural map
\[
k^\times/k^{\times p} \lang \prod_{\p\in  S} k_\p^\times/k_\p^{\times p}
\]
is injective, hence $V_S(k)=0$. Furthermore $V_\varnothing(k)$ is finite. Choosing to each nonzero element $\alpha$ of $V_\varnothing(k)$ a prime $\p\in S$ with $\alpha\notin k_\p^{\times p}$, we obtain a finite subset $T\subset S$ with $V_T(k)=0$.
\end{proof}

\begin{theorem}\label{infinitekpi1} Let $k$ be a number field and
let $S$ be a set of primes of $k$ of $p$-density $\varDelta^p(S)=1$.  Then $X\sm S$ is a $K(\pi,1)$ for $p$.
\end{theorem}

\begin{proof} Let $T\subset S$ be a finite subset. By \cite{NSW}, Thm.~9.2.2\,(ii), the natural map
\[
H^1_\et(X\sm (S\cup S_p)) \longrightarrow \prod_{\p\in T\cup S_p} H^1(k_\p)
\]
is surjective. A class in $H^1_\et(X\sm (S\cup S_p))$ which maps to zero in $H^1(k_\p)$ for all $\p \mid p$ is contained in $H^1_\et(X\sm S)$. Therefore, also the map
\[
H^1_\et(X\sm S) \longrightarrow \prod_{\p\in T} H^1(k_\p)
\]
is surjective.
Hence the maximal elementary abelian extension of $k$ in $k_S$ realizes the maximal elementary abelian extension of $k_\p$ in $k_\p(p)$ for all $\p\in S$. Applying the same argument to each finite subextension of $k$ in $k_S$, we conclude that $k_S$ realizes $k_\p(p)$ for all $\p\in S$. In particular,
\[
\prod_{\p \in S(k_S)} H^2((k_S)_\p) =0.
\]
Furthermore, by Proposition~\ref{sha2gleich0}, $\Sha^2(K,S(K))=0$ for all finite subextensions $K$ of $k$ in $k_S$. We obtain
\[
H^2_\et((X\sm S)_{k_S})=0.
\]
As there is no cohomology in dimension greater or equal~$3$, condition (iii) of Proposition~\ref{kpi1cond} is satisfied.
\end{proof}

In order to proceed, we make the following definitions.

\begin{definition} Let $S$ be a finite set of primes (of norm $\equiv 1 \bmod p$).
\begin{enumerate}
\item[\rm (i)] We say that $S$ is $p$-large if\/ $\Be_S(k,p)=0$.
\item[\rm (ii)] We put
\[
\delta^2_S(k)= h^2(X\sm S)- h^2(G_S)
\]
and call this number the $h^2$-defect of\/ $S$ (with respect to $p$).
\item[\rm (iii)] We denote by $k_S^{el}$ the maximal elementary abelian $p$-extension of $k$ inside~$k_S$.
\end{enumerate}
\end{definition}

\noindent
If $S$ is $p$-large, then $\Sha^2(k,S)=0$, and so, for any set $T\supset S$, the natural maps
$H^2(G_S)\to H^2(G_T)$ and $H^2_\et(X\sm S) \to H^2_\et(X\sm T)$ are injective.

\begin{lemma}\label{einsdazu}
Let $S$ be $p$-large and let $\p$ be a prime (of norm $\equiv 1 \bmod p$) which does not split completely in $k_S^{el}|k$. Then
\[
\delta^2_{S\cup \{\p\}}(k)\leq \delta^2_S(k).
\]
Furthermore, the natural map $H^2(G_{S\cup\{\p\}}) \mapr{} H^2(k_\p)$ is surjective.
\end{lemma}

\begin{proof}
Put $S'=S\cup \{\p\}$. By Theorem~\ref{globcoh}, the extension $k_{S'}^{el}|k$ is ramified at~$\p$. Therefore $k_{S'}^{el}$ realizes the maximal elementary abelian $p$-extension $k_\p^{el}$ of the local field $k_\p$, i.e.\ the map
\[
H^1(G_{S'}(p)) \longrightarrow H^1(k_\p)
\]
is surjective. As the cup-product $H^1(k_\p)\times H^1(k_\p) \to H^2(k_\p)$ is surjective, the natural map
\[
H^2(G_{S'}) \mapr{} H^2(k_\p)
\]
is also surjective. The statement of the lemma now follows from the commutative and exact diagram
\[
\renewcommand{\arraystretch}{1.5}
\begin{array}{ccccc}
0&\mapr{}& H^2(G_S)&\mapr{}& H^2_\et(X\sm S)\\
&&\injd{}&&\injd{}\\
0&\mapr{}& H^2(G_{S'})&\mapr{}& H^2_\et(X\sm S')\\
&&&&\surjd{}\\
&&&&H^2(k_\p).\\
\end{array}
\renewcommand{\arraystretch}{1}
\]
\end{proof}

\begin{lemma} \label{einsdazu1}
Let $S$ be $p$-large and let $\p$ be a prime. Let $T$ be a set of primes of $p$-density $\varDelta^p(T)=1$. Then there exists a prime $\p' \in T$ such that
\begin{enumerate}
\item[\rm (i)] $\p'$ does not split completely in $k_S^{el}|k$.
\item[\rm (ii)] $\p$ does not split completely in $k_{S\cup \{\p'\}}^{el}|k$.
\end{enumerate}
In particular, $\delta^2_{S\cup \{\p,\p'\}}(k)\leq \delta^2_S(k)$.
\end{lemma}

\begin{proof}
If $\p$ does not split completely in $k_S^{el}|k$, then condition (ii) is void. By assumption, $\varDelta^p(T)=\delta_{k(\mu_p)}(T(k(\mu_p))=1$. By \v Cebotarev's density theorem, we can find a prime $\P'\in T(k(\mu_p))$ which does not split completely in $k_S^{el}(\mu_p)|k(\mu_p)$. Then $\p'=\P'|_k$ satisfies~(i). Therefore we may assume that $\p$ splits completely in $k_S^{el}|k$. By class field theory, there exists an $s\in k^\times$ with

\smallskip
\begin{compactenum}
\item[\rm (a)] $v_\p(s)\equiv 1 \bmod p$,
\item[\rm (b)] $v_\q(s) \equiv 0 \bmod p$ for all $\q\notin S$, $\q \neq \p$, and
\item[\rm (c)] $s\in k_\q^{\times p}$ for all $\q \in S$.
\end{compactenum}

\smallskip\noindent
Since $S$ is $p$-large, $s$ is well-defined as an element in $k^\times/k^{\times p}$.  Now consider the extensions $k(\mu_p, \sqrt[p]{s})$ and $k_{S\cup\{\p\}}^{el}(\mu_p)$ of $k(\mu_p)$. The first one might be contained in the second (only if  $\zeta_p\in k$) but this does not matter.
Using \v Cebotarev's density theorem, we find $\P'\in T(k(\mu_p))$ such that $\Frob_{\P'}$ is non-zero in $Gal(k(\mu_p, \sqrt[p]{s})|k(\mu_p))$ and non-zero in $Gal(k_{S\cup\{\p\}}^{el}(\mu_p)|k(\mu_p))$. We put $\p'=\P'|_k$. Then $\p'$ does not split completely in $k_S^{el}|k$ and  $s\notin  k_{\p'}^{\times p}=k(\mu_p)_{\P'}^{\times p}$.

We claim that $\p$ does not split completely in $k_{S\cup \{\p'\}}^{el}|k$: Otherwise there would exist a $t \in k^\times$ satisfying conditions (a) -- (c) above and with $t \in k_{\p'}^{\times p}$. Since $s/t\in \Be_S(k)=0$, we obtain $s/t \in k^{\times p}$. Hence $s\in k_{\p'}^{\times p}$ giving a contradiction. Hence condition (i) and (ii) are satisfied.
\end{proof}

\begin{lemma}\label{extlemma1}
Let $S$ be a finite set of primes and let $T$ be a set of primes of $p$-density $\varDelta^p(T)=1$.  Then there exists a finite subset $T_1\subset T$  such that $S\cup T_1$ is $p$-large and such that the natural inclusion
\[
H^2(G_{S\cup T_1}(k)) \injr{} H^2_\et(X\sm (S\cup T_1))
\]
is an isomorphism.
\end{lemma}

\begin{proof}
We first move finitely many primes from $T$ to $S$ to achieve that $S$ is $p$-large. If $\delta^2_S(k)$ is zero, we are ready. Otherwise, consider the commutative diagram
\[
\renewcommand{\arraystretch}{1.5}
\begin{array}{ccc}
H^2(G_S)& \injr{}& H^2(G_{S\cup T})\\
\injd{}&&\mapd{\wr}\\
H^2_\et(X\sm S)&\injr{}& H^2_\et(X\sm (S\cup T))\\
\end{array}
\renewcommand{\arraystretch}{1}
\]
in which the right hand isomorphism follows from Theorem~\ref{infinitekpi1}. Let $x \in H^2_\et(X\sm S)$ but $x\notin H^2(G_S)$. Then there exists a finite subset $T_0 \subset T$ such that $x\in H^2(G_{S\cup T_0})$. Let $T_0=\{\p_1,\ldots,\p_n\}$. We choose $\p_1',\ldots,\p_n'\in T$ according to Lemma~\ref{einsdazu1} and put $T_1= \{\p_1,\ldots,\p_n, \p_1',\ldots,\p_n'\}$.
Then the natural map
\[
H^2(G_{S\cup T_1}) \mapr{\phi} \prod_{i=1}^n H^2(k_{\p_i}) \times \prod_{i=1}^n H^2(k_{\p_i'})
\]
is surjective. We have $H^2(G_S) \subset \ker(\phi)$ and also $x\in \ker(\phi)$. Hence
$\delta^2_{S\cup T_1}(k) < \delta^2_S(k)$. Iterating this process, we obtain a set with trivial $h^2$-defect.
\end{proof}

\section{\boldmath Review of mild pro-$p$-groups}

In the following we recall definitions and results from J.~Labute's paper \cite{La}. Only interested in our application, we are slightly less general than Labute.

\medskip
Let $R$ be a principal ideal domain and let $L$ be the free $R$-Lie algebra over $\xi_1,\ldots,\xi_n$, $n\geq 1$. We view $L$ as graded algebra where the degree of $\xi_i$ is $1$.
Let $\rho_1,\ldots,\rho_m$, $m\geq 1$, be homogeneous elements in $L$ with $\rho_i$ of degree $h_i$ and let $\r=(\rho_1,\ldots,\rho_m)$ be the ideal of $L$ generated by $\rho_1,\ldots,\rho_m$. Let $\g=L/\r$ and $U_\g$ be the universal enveloping algebra of $\g$.
Then $M=\r/[\r,\r]$ is a $U_\g$-module via the adjoint representation.

\begin{definition}
The sequence $\rho_1,\ldots,\rho_m$ is called {\bf strongly free} if $U_\g$ is a free $R$-module and $M=\r/[\r,\r]$ is the free $U_\g$-module on the images of $\rho_1,\ldots,\rho_m$ in~$M$.
\end{definition}

Let us consider the special case when $R=k[\pi]$ is the polynomial ring in one variable $\pi$ over a field $k$. Then $\bar L=L/\pi$ is a free $k$-Lie algebra and the next theorem shows that we can detect strong freeness by reduction. We denote the image in $\bar L$ of an element $\rho \in L$ by~$\bar \rho$.

\begin{theorem} {\rm(\cite{La}, Th.\,3.10)} \label{mildreduction}
The sequence $\rho_1,\ldots,\rho_m$ in $L$ is strongly free if and only if the sequence $\bar \rho_1,\ldots,\bar \rho_m$ is strongly free in $\bar L$.
\end{theorem}

Over fields, we have the following criterion for strong freeness. Let $R=k$ be a field, $X$  a finite set and $S\subset X$ a subset. Let $L(X)$ be the free Lie algebra over $X$ and let $\a$ be the ideal of $L(X)$ generated by the elements $\xi \in X\sm S$.
Put
\[
T=\{\ [\xi,\xi'] \ | \ \xi\in X\sm S,\ \xi' \in S  \} \subset \a.
\]

\begin{theorem}{\rm (\cite{La}, Th.\,3.3, Cor.\,3.5)} \label{mildoverfield} If $\rho_1,\ldots,\rho_m$ are homogeneous elements of $\a$ which lie in the $k$-span of\/ $T$ modulo $[\a,\a]$ and which are linearly independent over $k$ modulo $[\a,\a]$ then the sequence $\rho_1,\ldots,\rho_m$ is strongly free in $L$.
\end{theorem}

Let $p$ be an odd prime number and let $G$ be a pro-$p$-group. We consider the descending $p$-central series $(G_n)_{n\geq 1}$, which is defined recursively by
\[
G_1=G,\, G_{n+1}=G_n^p[G,G_n].
\]
The quotients $\gr_n(G)=G_n/G_{n+1}$, denoted additively, are $\F_p$-vector spaces. The graded vector space
\[
\gr(G)=\bigoplus_{n\geq 1} \gr_n(G)
\]
has a Lie algebra structure over the polynomial ring $\F_p[\pi]$, where multiplication by $\pi$ is induced by $x\mapsto x^p$ and the bracket operation for homogeneous elements is induced by the commutator operation in $G$, see \cite{NSW}, III,\,\S8. For $g\in G$, $g\neq 1$, let the natural number $h(g)$ be defined by
\[
g \in G_{h(G)},\ g\notin G_{h(G)+1}.
\]
This definition makes sense because $\bigcap_n G_n=\{1\}$, see \cite{NSW}, Prop.\,3.8.2. The image $\omega(g)$ of $g$ in $\gr_{h(g)}(G)$ is called the {\bf initial form} of~$g$.

\medskip
Let $F$ be the free pro-$p$-group over elements $x_1,\ldots,x_n$, $n\geq 1$. Then $h(x_i)=1$, $i=1,\ldots,n$, and
\[
L= \gr(F)
\]
is the free $\F_p[\pi]$-Lie algebra over $\xi_1,\ldots,\xi_n$, where $\xi_i=\omega(x_i)$, $i=1,\ldots,n$, see \cite{Lz}. Let $r_1,\ldots, r_m$, $m\geq 1$, be a sequence of elements in $F_2=F^p[F,F]\subset F$ and let $R=(r_1,\ldots,r_m)_F$ be the closed normal subgroup of $F$ generated by $r_1,\ldots,r_m$. Put $\rho_i=\omega(r_i) \in L$.

\begin{definition}
 A pro-$p$-group $G$ is called {\bf mild} if there exists a presentation
\[
1 \lang R \lang F \lang G \lang 1
\]
with $F$ a free pro-$p$-group on generators $x_1,\ldots,x_n$ and $R=(r_1,\ldots,r_m)_F$ such that the associated sequence  $\rho_1,\ldots,\rho_m$ is strongly free in $L=\gr(F)$.
\end{definition}

Essential for our application is the following property of mild pro-$p$-groups.

\begin{theorem}{\rm (\cite{La}, Th.\,1.2(c))} \label{mildcd2} \
If\/ $G$ is a mild pro-$p$-group, then $\cd G=2$.
\end{theorem}

Now let $G$ be a finitely presented pro-$p$-group and let
\[
1\lang R \lang F\lang G\lang 1
\]
be a minimal presentation, i.e.\ $F$ is the free pro-$p$-group on generators $x_1,\ldots, x_n$ where $n=\dim_{\F_p} H^1(G)$ and $R=(r_1,\ldots,r_m)_F$ with $m=\dim_{\F_p} H^2(G)$, cf.\ \cite{NSW},\,(3.9.5). Then the images $\xi_i=\omega(x_i)$, $i=1,\ldots, n$, of $x_1, \ldots, x_n$ are a basis of the $\F_p$-vector space  $F/F_2=H_1(F)=H_1(G)=G/G_2$.  For $y\in F_n$ and $a\in \Z_p$ the class of $y^a$ modulo $F_{n+1}$ only depends on the residue class $\bar a \in \F_p$ of~$a$.
Every element $r\in R\subset F_2$ has a representation
\[
r\equiv \prod_{j=1}^n (x_j^p)^{a_j} \  \cdot \prod_{1\leq k < l \leq n} [x_k,x_l]^{a_{kl}} \mod  F_3,
\]
where $a_j,a_{k,l} \in \F_p$. These coefficients are uniquely defined and can be calculated as follows. As $F$ is free, we have an isomorphism $H_2(G)=R_G^\ab/p$. Let $\bar r \in H_2(G)$ be the image of $r$ and
let $\chi_1,\ldots,\chi_n \in H^1(G)$ be the dual $\F_p$-basis of $\xi_1,\ldots,\xi_n$.

\begin{theorem}\label{coeff} \ $ a_{kl}=-\bar r(\chi_k\cup\chi_l)$ for $k<l$.
\end{theorem}

\noindent
For a proof see \cite{NSW}, Prop.\,3.9.13, which also gives a description of the $a_j$ using the Bockstein operator.

\bigskip
Using the results above, we obtain a criterion for mildness.

\begin{theorem} \label{mildkrit}
Let $G$ be a finitely presented pro-$p$-group. Assume there exists a basis $\chi_1,\ldots,\chi_n$ of $H^1(G)$, a basis $\bar r_1,\ldots,\bar r_m$ of $H_2(G)$ and a natural number $a$, $1\leq a < n$, such that the following conditions are satisfied
\begin{itemize}
\item[\rm (i)] $\bar r_i(\chi_k \cup \chi_l)=0$ for $a< k<l\leq n$ and $i=1,\ldots,m$.
\item[\rm (ii)] The $m{\times}a(n-a)$ matrix
\[
\Big(\bar r_i(\chi_k \cup \chi_l)\Big)_{i, (k,l)},\ 1\leq i\leq m,\ 1\leq k\leq a < l \leq n
\]
has rank $m$.
\end{itemize}
Then $G$ is a mild pro-$p$-group.
\end{theorem}

\begin{proof}
Let $\xi_1,\ldots,\xi_n\in H_1(G)$ be the dual basis of $\chi_1,\ldots, \chi_n$.
We choose a minimal presentation
\[
1\lang R \lang F\lang G\lang 1 \leqno (*)
\]
and generators $x_1,\ldots,x_n\in F$ mapping to $\xi_1,\ldots,\xi_n\in H_1(F)= H_1(G)$. Then we choose elements $r_1,\ldots,r_m\in R$ mapping to $\bar r_1,\ldots,\bar r_m\in R^\ab_G/p=H_2(G)$. Let $X= \{\xi_{1},\ldots,\xi_n\}$. Then  $L=\gr(F)$ is the free $\F_p[\pi]$-Lie algebra over $X$ and $\bar L=L/\pi$ is the free $\F_p$-Lie algebra over $X$. In order to show that $G$ is mild, we have to show that the initial forms $\rho_1,\ldots,\rho_m$ of $r_1,\ldots,r_m$ are a strongly free sequence in $L$. By Theorem~\ref{mildreduction} it suffices to show that
$\bar \rho_1,\ldots,\bar\rho_m \subset \bar L$ are a strongly free sequence.
By condition (ii) and Theorem~\ref{coeff}, we have $\bar \rho_1,\ldots,\bar\rho_m
 \in \gr_2(\bar L)=F_2/F_3F^p$.

Now we apply Theorem~\ref{mildoverfield} with $S=\{\xi_{a+1},\ldots,\xi_n\} \subset X$.
In the notation of this theorem, $\a$ is the ideal generated by $\xi_1,\ldots,\xi_a$ in $\bar L$ and
\[
T=\{[\xi_i,\xi_j]\ |\ 1\leq i\leq a,\ a+1\leq j \leq n\}.
\]
By condition (i) and Theorem~\ref{coeff}, we have $\bar \rho_i$ in the $\F_p$-span $H$ of $T$ modulo $[\a,\a]$. The elements of $T$ are a basis of $H$ and the coefficient matrix of $\bar \rho_1,\ldots,\bar \rho_m$ is up to sign the transpose of the matrix written in condition (ii). Hence $\bar \rho_1,\ldots,\bar \rho_m$ are linearly independent and, by Theorem~\ref{mildoverfield},  a strongly free sequence. This concludes the proof.
\end{proof}

\section{\boldmath Existence of $K(\pi,1)$'s}

Let $k$ be a number field and let $p$ be a prime number with $\mu_p\not \subset k$ and assume that $\Cl(k)(p)=0$. The exact sequence
\[
0\lang \O_k^\times \lang k^\times \mapr{(v_\q)_\q}\bigoplus_\q \Z \lang \Cl(k) \lang 0
\]
implies the exactness of
\[
0\lang \O_k^\times/p \lang k^\times/k^{\times p}\lang \bigoplus_\q \Z/p\Z \lang 0.
\]
Let $S=\{\p_1,\ldots,\p_m\}$ be a finite set of primes of norm $\equiv 1 \bmod p$.  We choose for $i=1,\ldots,m$ elements $s_i\in k^\times/k^{\times p}$ with $v_{\p_i}(s_i)\equiv 1 \bmod p$ and $v_\q(s_i)\equiv 0 \bmod p$ for all primes $\q\neq \p_i$ of $k$. Let furthermore, $e_1,\ldots,e_r$, $r=r_1+r_2-1$, be a basis of $\O_k^\times/p$.

\medskip
Consider the field
\[
K= k(\mu_p, \sqrt[p]{s_1},  \ldots, \sqrt[p]{s_m},  \sqrt[p]{e_1},  \ldots, \sqrt[p]{e_r}).
\]
An inspection of the ramification behaviour shows that $Gal(K|k(\mu_p))$ has the Galois group
$(\Z/p\Z)^{m+r}$: Indeed, $k(\mu_p,\sqrt[p]{e_1}, \ldots, \sqrt[p]{e_r})|k(\mu_p)$  is unramified outside $S_p$ and has Galois group $(\Z/p\Z)^{r}$ by Kummer theory. Adjoining $\sqrt[p]{s_i}$, $i=1,\ldots,m$, yields a cyclic extension of degree $p$ which is unramified outside $S_p\cup \{\p_i\}$ and ramified at $\p_i$.

Since $\mu_p\not\subset k$, the extensions $k_S^{el}(\mu_p)|k(\mu_p)$ and $K|k(\mu_p)$ lie in different eigenspaces for the action of $Gal(k(\mu_p)|k)$. Therefore  $Kk_S^{el}|k(\mu_p)$ has Galois group $(\Z/p\Z)^{m+r+n}$, with
$n=\dim_{\F_p} Gal(k_S^{el}|k)=\dim_{\F_p}H^1(G_S)$.

\medskip
Assume now that we are given

\medskip\noindent
\begin{compactitem}
\item a set of primes $T$ of $k$ with $T \cap S =\varnothing$ and with $p$-density $\varDelta_k^p(T)=1$,
\item a nonzero element
$F\in Gal(k_S^{el}|k)=Gal(k_S^{el}(\mu_p)|k(\mu_p))$,

\item to each $\p_i$, $i=1,\ldots,m$, a condition $C_i$ which says ``split'' or ``inert''.
\end{compactitem}

\medskip\noindent
By \v Cebotarev's density theorem applied to the extension $Kk_S^{el}(\mu_p)|k(\mu_p)$, we find a prime $\P \in T(Kk_S^{el}(\mu_p))$ such that

\medskip\noindent
\begin{compactitem}
\item the image of $\Frob_\P$ in $Gal({k(\mu_p,   \sqrt[p]{e_1}, \ldots, \sqrt[p]{e_r})|k(\mu_p)})$ is trivial,
\item the image of $\Frob_\P$ in $Gal({k(\mu_p, \sqrt[p]{s_i})|k(\mu_p)})$ is trivial if $C_i$ is ``split'' and nontrivial otherwise, and
\item the image of $\Frob_\P$ in $Gal(k_S^{el}(\mu_p)|k(\mu_p))$ is equal to $F$.
\end{compactitem}

\medskip\noindent
Let $\p\in T$ be the restriction of $\P$ to $k$. Then the natural map
$\O_k^\times/p \to k_\p^\times/k_\p^{\times p}$ is the zero map. Since $_p\Cl(k)=0$,  we obtain $\O_k^\times/p \stackrel{\sim}{\to} V_\varnothing(k)=V_{\{\p\}}(k)$. By Theorem~\ref{globcoh}, $k_{\{\p\}}^{el}|k$ is cyclic of order $p$ and $\p$ is ramified in this extension.
Recall that $H^1_\nr(G_\p)$ is defined as the exact annihilator of the inertia group $T_\p(k_\p^{el}|k_\p)\subset H_1(G_\p)$ in the natural pairing
\[
H_1(G_\p) \times H^1(G_\p) \longrightarrow \F_p.
\]
Dually, $T_\p(k_\p^{el}|k_\p)$ is the exact annihilator of $H^1_\nr(G_\p)$. The equation $T_\p(k_{\{\p\}}^{el}|k)=Gal(k_{\{\p\}}^{el}|k)$
yields an isomorphism
\[
 H^1(G_{\{\p\}}) \mapr{\sim} H^1(G_\p)/H^1_\nr(G_\p).
\]
By class field theory, $\p_i$ splits in $k_{\{\p\}}^{el}|k$ if and only if there exists an element $s'_i\in k^\times/k^{\times p}$ with $v_{\p_i}(s'_i)\equiv 1 \bmod p$, $v_\q(s'_i)\equiv 0 \bmod p$ for all $\q\neq \p_i$ and $s'_i\in k_\p^{\times p}$.
Then $s'_i/s_i$ lies in $\O_k^\times/p$, and therefore $s_i\in k_\p^{\times p}$.
Hence $\p_i$ splits in $k_{\{\p\}}^{el}|k$ if and only if $s_i$ is a $p$-th power in $k_\p$.
On the other hand,  by our choice of $\P$,  $s_i$ is a $p$-th power in $k_\p$ if and only if $C_i$ is ``split''. Therefore the following holds:

\medskip\noindent
\begin{compactitem}
\item the natural map $\O_k^\times/p\to k_\p^\times/k_\p^{\times p}$ is the zero map,
\item $\Frob_\p=F \in Gal(k_S^{el}|k)$,
\item $k_{\{\p\}}^{el}|k$ is cyclic of order $p$,
\item  each $\p_i$, $i=1,\ldots,m$, satisfies condition $C_i$ in $k_{\{\p\}}^{el}|k$.
\end{compactitem}

\medskip\noindent
Now assume that $\Be_{S\sm \{\q\}}(k)=0$ for all $\q\in S$, in particular, $S$ is $p$-large.
Then all $\p_i\in S$ ramify in $k_S^{el}|k$ and the $1$-dimensional subspaces $T_{\p_i}(k_S^{el}|k)$, $i=1,\ldots,m$, in $H_1(G_S)$ are pairwise different and generate $H_1(G_S)$.
Furthermore assume that  $\delta_S^2(k)=0$. As $\p$ does not split completely in $k_S^{el}|k$,  Lemma~\ref{einsdazu} implies $\delta_{S\cup \{\p\}}^2(k)=0$.
Since $\mu_p\not\subset k$ and by Theorem~\ref{globcoh}, the natural maps $H^2(G_S) \to \prod_{\q\in S} H^2(G_\q)$ and $H^2(G_{S\cup \{\p\}}) \to \prod_{\q\in S\cup \{\p\}} H^2(G_\q)$ are isomorphisms. We denote the $\q$-component of a global cohomology class $\alpha$ by $\alpha_\q$.

\medskip\noindent
Next we fix a primitive $p$-th root of unity in $k(\mu_p)$ and to each $\p_1,\ldots,\p_m$ a prolongation to $k(\mu_p)$. After this choice we have identifications $\mu_p(k_\p)=\mu_p((Kk_S^{el})_\P)\cong \F_p$ and $\mu_p(k_{\p_i})\cong \F_p$, $i=1,\ldots, m$.
In particular, we have an isomorphism $H^2(G_\p)=H^2(G_\p, \mu_p)=\F_p$, and similarly for the $\p_i$. Via these isomorphisms we consider the $\q$-component $\alpha_\q$ of a class $\alpha\in H^2(G_{S\cup\{\p\}})$ as an element in $\F_p$.
Let
\[
\pi_\p \in H^1(G_\p)/H^1_\nr(G_\p)= H^1(G_\p,\mu_p)/H^1_\nr(G_\p,\mu_p)= k_\p^\times/U_\p k_\p^{\times p}
\]
be the image of a uniformizer and let $\chi_\p\in H^1(G_{\{\p\}})$ be the unique pre-image. We denote the image of $\chi_\p$ in $H^1(G_{S\cup \{\p\}})$ by the same letter. Thus $\chi_\p$ maps to $\pi_\p$ under the natural map $H^1(G_{S\cup \{\p\}}) \to H^1(G_\p)/H^1_\nr(G_\p)$.
Consider the exact pairing
\[
H^1_\nr(G_\p) \times H^1(G_\p)/H^1_\nr(G_\p) \to H^2(G_\p)=\F_p,
\]
which is induced by local Tate duality, see \cite{NSW}, Thm.\,7.2.15. Let $\delta: k_\p^\times/k_\p^{\times p} \stackrel{\sim}{\to} H^1(G_\p)$ be the boundary isomorphism of the Kummer sequence and let $\mathit{rec}: k_\p^\times/k_\p^{\times p} \stackrel{\sim}{\to} H_1(G_\p)$ be the mod-$p$ reciprocity map. Put $\varphi= \mathit{rec}\circ \delta^{-1}$. Then the image of $\chi_\p$ under the composition
\[
H^1(G_S) \mapr{} H^1(G_\p) \mapr{\phi} H_1(G_\p) \mapr{} H_1(G_S)
\]
is $\Frob_\p$, the Frobenius automorphism of the unramified prime~$\p$ in $k_S^{el}|k$.
By  \cite{NSW}, Prop.\,7.2.13\ \footnote{This proposition contains a sign error, see the errata file on the author's homepage.}), the diagram
\[
\begin{array}{ccccc}
H^1(G_\p)&\times & H^1(G_\p)&\mapr{\cup}&H^2(G_\p)\\
\eqd&&\mapd{\hspace{-7pt} \wr \ \varphi}&&\mapd{\hspace{-7pt}\wr \ \mathit{inv}}\\
H^1(G_\p)&\times & H_1(G_\p)&\mapr{\mathit{can}}&\F_p
\end{array}
\]
commutes.
We obtain for any $\chi\in H^1(G_S)\subset H^1(G_{S\cup \{\p\}})$ the following formula for the $\p$-component of $\chi\cup \chi_\p \in H^2(G_{S\cup\{\p\}})$:
\[
(\chi\cup \chi_\p)_\p= \chi(\Frob_\p).
\]
The image of $\chi_\p$ in $H^1 (G_{\p_i})$ obviously lies in the subgroup $H^1_\nr (G_{\p_i})$. By the same argument, noting that the cup-product is anti-symmetric, we obtain the equality
\[
(\chi \cup \chi_\p)_{\p_i}= -\chi_\p(\Frob_{\p_i}),
\]
for any $\chi\in H^1(G_S)$ mapping to $\pi_{\p_i}\in H^1(G_{\p_i})/H^1_\nr(G_{\p_i})$, where $\Frob_{\p_i}$ is the Frobenius automorphism of the unramified prime $\p_i$ in $k_{\{\p\}}^{el}|k$. As $\chi_\p$ is unramified at $\p_i$, the element $(\chi \cup \chi_\p)_{\p_i}$ depends only on the image of $\chi$ in the one-dimensional $\F_p$-vector space $H^1(G_{\p_i})/H^1_{nr}(G_{\p_i})$.  Since  $\p_i$ ramifies in $k_S^{el}|k$, the map $H^1(G_S)\to \F_p$,  $\chi \mapsto (\chi \cup \chi_\p)_{\p_i}$ is the linear form associated to an element $t_i \in T_{\p_i}(k_S^{el}|k)\subset H_1(G_S)$.

\medskip
Summing up and using the notation and choices above, we obtain the
\begin{lemma} \label{search}
Let $k$ be a number field and let $p$ be a prime number with $\mu_p\not\subset k$ and $\Cl(k)(p)=0$. Let $S=\{\p_1,\ldots, \p_m\}$ be a finite $p$-large set of primes and assume $\delta^2_S(k)=0$. Let for $i=1,\ldots,m$ elements $\varepsilon_i\in \{0,1\}$ and for $i=1,\ldots,n$ elements $d_i \in \F_p$ be given, where not all $d_i$ are zero. Let $\chi_1,\ldots,\chi_n$ be a basis of $H^1(G_S)$. Furthermore, let $T$ be a set of primes of $p$-density $\varDelta_p(T)=1$ and with $T\cap S=\varnothing$.

Then there exists a prime $\p\in T$  such that the following conditions hold with respect to the identifications $H^2(G_{\p_i})=\F_p$,  $i=1,\ldots,m$,  and $H^2(G_\p)=\F_p$.

\medskip\noindent
\begin{compactitem}
\item $\p$ does not split completely in $k_S^{el}|k$,
\item $k_{\{\p\}}^{el}|k$ is cyclic of order $p$,
\item $\chi_1,\ldots, \chi_n,\chi_\p$ is a basis of $H^1(G_{S\cup\{\p\}})$,
\item $(\chi_i \cup \chi_\p)_\p=d_i$ for  $i=1,\ldots,n$,
\item For  $i=1,\ldots,m$ we have $c_i=0$ if and only if  $\varepsilon_i=0$, where $c_i\in T_{\p_i}(k_S^{el}|k)\subset H_1(G_S)$ represents the map $H^1(G_S) \to \F_p$, $\chi \mapsto (\chi \cup \chi_\p)_{\p_i}$.
\end{compactitem}

\end{lemma}

Now we are able to prove the following result, which is unessentially sharper than Theorem~\ref{exthmi} of the introduction.

\begin{theorem}\label{existenceofkpi1}
Let $k$ be a number field and let $p$ be a prime number such that
\[
\mu_p\not\subset k\ \text{ and }\ \Cl(k)(p)=0.
\]
Let $S$ be a finite set of primes of\/ $k$ and let $T$ be a set of primes of $p$-density $\varDelta^p(T)=1$. Then there exists a finite subset $T_1\subset T$ such that $\Spec(\O_k)\sm (S\cup T_1)$ is a $K(\pi,1)$ for\/ $p$.
\end{theorem}

\begin{proof}
We may suppose that $T\cap S=\varnothing$.
After moving finitely many primes of $T$ to $S$, we may assume that the following conditions hold:

\medskip\noindent
\begin{compactitem}
\item $\Be_{S\sm \{\p\}}(k)=0$ for all $\p\in S$,\smallskip
\item $\delta_S^2(k)=0$.
\end{compactitem}

\medskip\noindent
Now let $S=\{\p_1,\ldots,\p_m\}$. Then $m=h^2(G_S)$. Let $n=m-r=h^1(G_S)$. We will achieve the $K(\pi,1)$-situation by adding $m$ further primes to $S$.

We choose any basis $\chi_1,\ldots, \chi_n$ of $H^1(G_S)$. Let $t_1,\ldots,t_m$ be generators of the inertia groups $T_{\p_i}(k_S^{el}|k)\subset H_1(G_S)$.  Now we add a prime $\p_{m+1}$ in the following way:

Let $i_1 \in \{1,\ldots, n\}$ be an index such that $\chi_{i_1}(t_1)\neq 0$, and let $i_1' \in \{1,\ldots, n\}$, $i_1'\neq i_1$, be any other index. Now, according to Lemma~\ref{search},  we put the conditions
\[
\varepsilon_1=1  \text{ and } \varepsilon_i=0 \text{ for } i\in \{2,\ldots,m\},
\]
\[
 d_{i_1'}=1 \text{ and } d_i=0 \text{ for } i\in \{1,\ldots,n\},\ i\neq i_1'
\]
to choose a prime $\p_{m+1}\in T$  such that for $i=1,\ldots, n$
\[
(\chi_i \cup \chi_{\p_{m+1}})_{\p_1}= \lambda_1 \chi_i(t_1), \lambda_1 \in \F_p^\times,\  (\chi_i \cup \chi_{\p_{m+1}})_{\p_j}=0,\ j=2,\ldots,m
\]
\[
\text{ and } (\chi_i \cup \chi_{\p_{m+1}})_{\p_{m+1}}=d_i.
\]
Then in the matrix
\[
\Big( (\chi_i\cup \chi_{\p_{m+1}})_{\p_j}\Big)_{\substack{
i=1,\ldots,n\\ j=1,\ldots, m+1}}
\]
the $i_1$-line has entry $\neq 0$ at $(i_1,1)$ and all other entries zero, while the $i_1'$-line has some entry at $(i_1',1)$, the entry $1$ at $(i_1',m+1)$ and all other entries zero.

In order to proceed, we put $\chi_{n+1}=\chi_{\p_{m+1}}$ and choose an index $i_2\in \{1,\ldots,n\}$ with $\chi_{i_2}(t_2)\neq 0$ and any  $i'_2\in \{1,\ldots,n\}$ with $i'_2\neq i_2$. We choose conditions as before, completed by $\varepsilon_{m+1}=0$ and $d_{n+1}=0$. Then we choose $\p_{m+2}$ according to Lemma~\ref{search} and such that in the matrix
\[
\Big( (\chi_i\cup \chi_{\p_{m+2}})_{\p_j}\Big)_{\substack{
i=1,\ldots,n\\ j=1,\ldots, m+2}}
\]
the $i_2$-line has entry $\neq 0$ at $(i_2,2)$ and all other entries zero, while the $i_2'$-line has some entry at $(i_2',2)$, the entry $1$ at $(i_2',m+2)$ and all other entries zero. In addition, our choice implies
\[
(\chi_{\p_{m+1}} \cup \chi_{\p_{m+2}})_{\p_{m+1}}=0= (\chi_{\p_{m+1}} \cup \chi_{\p_{m+2}})_{\p_{m+2}}
\]
As $\chi_{\p_{m+1}}$ and $\chi_{\p_{m+2}}$ are unramified at $\p_1,\ldots,\p_m$ by construction, we have furthermore $(\chi_{\p_{m+1}} \cup \chi_{\p_{m+2}})_{\p_{i}}=0$ for $i=1,\ldots,m$.

\medskip
Now we proceed to construct $\p_{m+3},\ldots, \p_{2m}$ in a similar way, and apply Theorem~\ref{mildkrit} with $a=m$. For each $j$, the $j$-th of the $m$-steps in the construction produced the two lines $((i_j,j), - )$ and $((i'_j,j), - )$ in the $nm{\times}2m$-matrix
\[
\Big( (\chi_i\cup \chi_{\p_{j}})_{\p_k}\Big)_{\substack{
i=1,\ldots,n, j=m+1,\ldots, 2m\\
k=1,\ldots, 2m}}
\]
According to our choices these $2m$ lines are linearly independent, hence the matrix has rank $2m$. Putting $T_1=\{\p_{m+1},\ldots,\p_{2m}\}$, we conclude by Theorem~\ref{mildkrit} that $G_{S\cup T_1}$ is a mild pro-$p$-group. Hence $\cd G_{S\cup T_1}=2$ by Theorem~\ref{mildcd2}. By Lemma~\ref{einsdazu}, we didn't produce new $h^2$-defect during our construction, hence $\delta^2_{S\cup T_1}(k)=0$. As the \'{e}tale cohomology is trivial in dimension $\geq 3$,  we conclude that the homomorphisms
\[
\phi_i: H^i(G_{S\cup T_1},\Z/p\Z) \lang H^i_et(\Spec(\O_k)\sm (S\cup T_1),\Z/p\Z)
\]
are isomorphisms for all $i\geq 0$. Hence condition (v) of Proposition~\ref{kpi1cond} is satisfied.
\end{proof}

\section{\boldmath Consequences of the $K(\pi,1)$-property}

In this section we assume that $S$ is finite and we exclude the case $S=\varnothing$ from our considerations. Keeping all conventions made before, we assume

\medskip
\begin{minipage}{11cm}{\it  $p\neq 2$ or $k$ is totally imaginary and
 $S$ is a non-empty finite set of nonarchimedean primes $\p$ with norm \hbox{$N(\p)\equiv 1 \bmod p$.}}
\end{minipage}

\begin{lemma} \label{fab} $G_S$ is a fab-group, i.e.\
the abelianization $U^\ab$ of every open subgroup $U$ of $G_S$ is finite.
\end{lemma}

\begin{proof} Let $U \subset G_S$ be an open subgroup. The abelianization $U^{\ab}$ of $U$ is a finitely generated abelian pro-$p$-group. If $U^{\ab}$ were infinite, it would have a quotient isomorphic to $\Z_p$, which by Galois theory corresponds to a $\Z_p$-extension $K_\infty$ of the number field $K=k_S^U$ inside $k_S$.
By \cite{NSW}, Theorem 10.3.20 (ii), a $\Z_p$-extension of a number field is ramified at at least one prime dividing $p$. This contradicts $K_\infty\subset k_S$ and we conclude that $U^{\ab}$ is finite.
\end{proof}

The group theoretical structure of the local Galois groups is well known.

\begin{proposition}\label{localstrukt}
Let $\p\in S$. Then $\Gal(k_\p(p)|k_\p)$ is the pro-$p$-group on two generators $\sigma, \tau$ subject to the relation $\sigma\tau\sigma^{-1}=\tau^q$. The element $\tau$ is a generator of the inertia group, $\sigma$ is a Frobenius lift and $q=N(\p)$.
\end{proposition}

\begin{proof}
This follows from \cite{NSW}, Thm.~7.5.2 by passing to the maximal pro-$p$-factor group.
\end{proof}

We obtain the following corollary.

\begin{corollary}\label{infinite_index}
Assume that $G_S$ is infinite. Then, for each $\p\in S$, the decomposition group $G_\p$ of\/ $\p$ in $G_S$ has infinite index.
\end{corollary}

\begin{proof} The decomposition group $G_\p$ is a quotient of the local Galois group $\Gal(k_\p(p)|k_\p)$.  If $G_\p \subset G_S$ would have finite index, it would be an infinite fab-group by Lemma~\ref{fab}.  By Proposition~\ref{localstrukt}, each infinite quotient of $\Gal(k_\p(p)|k_\p)$ has a surjection to $\Z_p$ and is therefore not a fab-group. This contradiction shows that $G_\p$ has infinite index in $G_S$.
\end{proof}

\bigskip
The next proposition classifies the degenerate $K(\pi,1)$-case.

\begin{proposition}\label{degen}
$X\sm S$ is a $K(\pi,1)$  and  $G_S=1$ if and only if  $S=\{\p\}$ consists of a single prime and one of the following cases occurs.
\begin{itemize}
\item[\rm (a)] $p=2$, $k\neq\Q(\sqrt{-1})$ is imaginary quadratic, $2\nmid h_k$ and \hbox{$N(\p)\not \equiv 1 \bmod 4$},
\item[\rm (b)] $p=2$, $k=\Q(\sqrt{-1})$ and $N(\p)\not \equiv 1 \bmod 8$,
\item[\rm (c)] $p=3$, $k=(\Q\sqrt{-3})$ and $N(\p) \not \equiv 1 \bmod 9$.
\end{itemize}
\end{proposition}

\begin{proof}  Assume $G_S=1$ and that $X \sm S$ is a $K(\pi,1)$. Then $H^i_\et(X\sm S)=0$ for all $i\geq 1$. In particular, $p\nmid h_k$. By Theorem~\ref{globcoh}, $h^2(X\sm S)=0$ implies $\delta=1$, $\# S=1$ and $V_S=0$. Then, using $h^1(X\sm S)=0$, we obtain $r=1$. As $\delta=1$, the following possibilities remain

\begin{itemize}
\item[\rm (a)] $p=2$, $k\neq\Q(\sqrt{-1})$ is imaginary quadratic and $2\nmid h_k$,
\item[\rm (b)] $p=2$, $k=\Q(\sqrt{-1})$,
\item[\rm (c)] $p=3$, $k=(\Q\sqrt{-3})$.
\end{itemize}
In all cases, Proposition~\ref{VSchange} yields an isomorphism
$\O_k^\times/p \stackrel{\sim}{\to} V_\varnothing$.  The second exact sequence of Proposition~\ref{VSchange} and the isomorphism $U_\p/p \cong U_\p k_\p^{\times p}/k_\p^{\times p}$ imply
\[
0=V_S=\ker \big( \O_k^\times/p \to U_\p/p  \big).
\]
Note that $\O_k^\times/p$\/ is one-dimensional.
In case (a), the unit $-1$ is a generator of $\O_k^\times/2$ which must not be a square in  $U_\p$, implying $N(\p)\not\equiv 1 \bmod 4$.   In case (b), $\sqrt{-1}$ is a generator, and in case (c), a generator is given by $\zeta_3=\frac{1}{2}(-1+\sqrt{-3})$. The assertions  in the cases (b) and (c) follow similarly. Conversely, assume we are in case (a), (b) or (c). Then we can reverse the given arguments and obtain $h^i(X\sm S)=0$ for all $i\geq 1$.
\end{proof}

\begin{theorem}
Assume  $G_S\neq 1$ and that $X\sm S$ is a $K(\pi,1)$. Then the following holds.

\begin{itemize}
\item[\rm (i)] $\cd G_S=2$, $\scd G_S=3$.
\item[\rm (ii)] $G_S$ is a duality group (of dimension~$2$).
\end{itemize}
\end{theorem}

\begin{proof}
By Lemma~\ref{fab} and Corollary \ref{condforkpi1}, $G_S$ is a fab-group and $\cd\,G_S\leq 2$. Now the assertions follow in a purely group-theoretical way:

As $G_S\neq 1$ and $G_S^{\ab}$ is finite, $G_S$ is not free, and we obtain $\cd G_S=2$. By \cite{NSW}, Proposition 3.3.3, it follows that $\scd G_S \in \{2,3\}$. Assume  $\scd G=2$. We consider the $G_S$-module
\[
D_2(\Z):=\varinjlim_{U} U^{\ab},
\]
where the limit runs over all open normal subgroups $U\lhd G_S$ and for $V \subset U$ the transition map is the transfer $\text{Ver}\colon U^{\ab} \to V^{\ab}$, i.e.\ the dual of the corestriction map $\text{cor}\colon H^2(V,\Z) \to H^2(U,\Z)$ (see \cite{NSW}, I, \S5).
By \cite{NSW}, Theorem 3.6.4 (iv), we obtain $G_S^{\ab}=D_2(\Z)^{G_S}$. On the other hand, $U^{\ab}$ is finite for all $U$ and the group theoretical version of the Principal Ideal Theorem (see \cite{Ne}, VI, Theorem 7.6) implies $D_2(\Z)=0$. Hence $G_S^{\ab}=0$ which implies $G_S=1$ producing a contradiction. Hence $\scd G_S=3$.

It remains to show that $G_S$ is a duality group. By \cite{NSW}, Theorem 3.4.6, it suffices to show that the terms
\[
D_i(G_S,\Z/p\Z):= \varinjlim_{U} H^i(U,\Z/p\Z)^\vee
\]
are zero for $i=0,1$. Here $U$ runs through the open subgroups of $G_S$, $\scriptstyle \vee$ denotes the Pontryagin dual and the transition maps are the duals of the corestriction maps. For $i=0$, and $V\subsetneqq U$, the transition map
\[
\text{cor}^\vee\colon \Z/p\Z=H^0(V,\Z/p\Z)^\vee \to H^0(U,\Z/p\Z)^\vee=\Z/p\Z
\]
is multiplication by $(U:V)$, hence zero. Therefore $D_0(G_S,\Z/p\Z)=0$, as $G_S$ is infinite. Furthermore,
\[
D_1(G_S,\Z/p\Z)=\varinjlim_U U^{\ab}/p=0
\]
by the Principal Ideal Theorem. This finishes the proof.
\end{proof}

In order to proceed, we introduce some notation in order to deal with the case of infinite extensions.  For a (possibly infinite) algebraic extension $K$ of $k$ we denote by $S(K)$ the set of prolongations of primes in $S$ to $K$. The set $S(K)$ carries a profinite topology in a natural way. Now assume that $M|K|k$ is a tower of Galois extensions.  We denote the inertia group of a prime ${\mathfrak p}\in S(K)$ in the extension $M|K$ by $T_{\mathfrak p}(M|K)$. For $i\geq 0$ we write
\[
\ressum_{{\mathfrak p} \in S(K)} H^i(T_{\mathfrak p}(M|K), \Z/p\Z)
\stackrel{df}{=}
\varinjlim_{k'\subset K} \bigoplus_{{\mathfrak p} \in S(k')} H^i(T_{\mathfrak p}(M|k'), \Z/p\Z),
\]
where the limit on the right hand side runs through all finite subextensions $k'$ of $k$ in $K$. The $\Gal(K|k)$-module $\ressumsmall_{{\mathfrak p} \in S(K)} H^i(T_{\mathfrak p}(M|K), \Z/p\Z)$ is the maximal discrete submodule of the  product $\prod_{{\mathfrak p} \in S(K)} H^i(T_{\mathfrak p}(M|K), \Z/p\Z) $.

\medskip
Whenever we deal with local terms associated to the elements of $S(K)$ (e.g.\ \'{e}tale cohomology groups) we use restricted sums, which are, in the same manner as above, defined as the inductive limit over the similar terms associated to all finite subextensions of $k$ in $K$.

\bigskip
A natural question is how far we get locally at the primes in $S$ when going up to~$k_S$.
\begin{proposition} \label{wieweitlokal}
Assume that $X\sm S$ is a $K(\pi,1)$  and that $G_S\neq 1$. Then $k_S$ realizes the maximal unramified $p$-extension of $k_\p$ for all $\p\in S$, i.e.\
\[
k_\p^{nr}(p) \subset (k_S)_\p \quad \text{ for all } \p\in S.
\]
If\/ $\p \in S$ ramifies in $k_S$, then $(k_S)_\p=k_\p(p)$, i.e.\ $k_S$ realizes the maximal $p$-extension of\/ $k_\p$.
\end{proposition}

\begin{proof} For an integral normal scheme $Y$ we write $Y_L$ for the normalization of $Y$ in an algebraic extension $L$ of its function field.
Then $(X \sm S)_{k_S}$ is the universal pro-$p$ covering of $X\sm S$.
We consider the following part of the excision sequence for the pair $(X_{k_S},  (X \sm S)_{k_S}) $
\[
H^2_\et((X\sm S)_{k_S}) \to \ressum_{\p \in S(k_S)} H^3_\p((X_{k_S})_\p) \to H^3_\et(X_{k_S}) .
\]
As $G_S$ is infinite, Lemma~\ref{h3oben} implies $H^3_\et(X_{k_S})=0$. By condition (iii) of Proposition~\ref{kpi1cond} we have $H^2_\et((X\sm S)_{k_S})=0$.  Hence $H^3_\p((X_{k_S})_\p)=0$ for all $\p \in S(k_S)$.
As $H^i_\et ((X_{k_S})_\p)=0$ for $i\geq 2$, we obtain
\[
H^3_\p((X_{k_S})_\p)\cong H^2((k_S)_\p),
\]
where the group on the right hand side is Galois cohomology with values in $\Z/p\Z$.
As $\mu_p\subset k_\p$ by assumption, the vanishing of $H^2((k_S)_\p)$ implies $p^\infty\mid [(k_S)_\p:k_\p]$. In other words, the decomposition group $G_\p(k_S|k)$ of each $\p\in S$ is infinite. As a subgroup of $G_S$, it has cohomological dimension $\leq 2$. Furthermore, $G_\p(k_S|k)$ is a factor group of the local Galois group $\Gal(k_\p(p)|k_\p)$, which, by Proposition~\ref{localstrukt}, has only three quotients of cohomological dimension less or equal to $2$: itself, the trivial group and the Galois group of the maximal unramified $p$-extension of $k_\p$. Hence $k_\p^{nr}(p) \subseteq (k_S)_\p$ and $(k_S)_\p=k_\p(p)$ if $\p$ ramifies in $k_S$.
\end{proof}

\bigskip
In order to deduce Theorem~\ref{1.2}, it remains to show that each $\p\in S$ ramifies in $k_S$. The following lemma provides a first step.

\begin{lemma} \label{whenramifies} Let $\p \in S$ be a prime and let $S'=S\sm \{\p\}$. Assume that the natural injection $V_{S} \hookrightarrow V_{S'}$ is an isomorphism. Then $\p$ ramifies in $k_S$.
\end{lemma}

\begin{proof}
Since the map $H^1(G_S)\to H^1_\et(X\sm S)$ is an isomorphism,  Theorem~\ref{globcoh} implies
\[
\dim_{\F_p}H^1(G_S)= 1 + \# S - \delta + \dim_{\F_p} V_S - r,
\]
and the same formula holds with $S$ replaced by $S'$. Hence
\[
\dim_{\F_p}H^1(G_S) = \dim_{\F_p}H^1(G_{S'}) +1\,.
\]
In particular, $G_{S'}$ is a proper quotient of $G_S$ and therefore $\p$ ramifies in $k_S$.
\end{proof}

\begin{corollary}
Assume that $X\sm S$ is a $K(\pi,1)$  and that $G_S\neq 1$. Let $\p \in S$ be a prime and let $S'=S\sm \{\p\}$. Assume that $V_{S'}=0$. Then $(k_S)_\p=k_\p(p)$.
\end{corollary}

\noindent
{\bf Remark:} If $V_\varnothing =0$, then the given criterion applies to any set $S$ and each $\p \in S$. This was used in \cite{S} for $k=\Q$ and in \cite{Vo} for imaginary quadratic number fields. If the unit rank of $k$ is non-zero, then $V_\varnothing\neq 0$ and the criterion applies only to sufficiently large sets $S$.

\section{Enlarging the set of primes}
Next we consider the problem of enlarging the set $S$.

\begin{proposition} \label{enlarge}
Let $S\subset S'$ be finite sets of primes of norm congruent to $1$ modulo~$p$.  Assume that $X\sm S$ is a $K(\pi,1)$ and that $G_S\neq 1$. Further assume that each $\q\in S'\sm S$ does not split completely in $k_S$. Then the following holds.

\begin{itemize}
\item[\rm (i)] $X\sm S'$ is a $K(\pi,1)$.
\item[\rm (ii)] $(k_{S'})_\q=k_\q(p)$ for all $\q \in S'\sm S$.
\end{itemize}
Furthermore,
$
H^i(\Gal(k_{S'}|k_S))=0
$
for $i\geq 2$. For $i=1$ we have a natural isomorphism
\[
H^1(\Gal(k_{S'}|k_S)) \cong \ressum_{{\mathfrak p} \in S'\sm S (k_S))} \!\!\!\! H^1(T_{\mathfrak p}(k_{S'}|k_S)),
\]
In particular, $\Gal(k_{S'}|k_S)$ is a free pro-$p$-group.
\end{proposition}

 \begin{proof} Let $\q \in S'\sm S$. Since $\q$ does not split completely in $k_S$ and since $\cd G_S=2$, the decomposition group of $\q$ in $k_S|k$ is a non-trivial and torsion-free quotient of $\Z_p\cong G(k_\q^{nr}(p)|k_\q)$. Therefore $(k_S)_\q$ is the maximal unramified $p$-extension of $k_\q$. We denote the normalization of an integral normal scheme $Y$ in an algebraic extension $L$ of its function field by $Y_L$.
 Then $(X\sm S)_{k_S}$ is the universal pro-$p$ covering of $X\sm S$.  We consider the \'{e}tale excision sequence for
the pair $((X \sm S)_{k_S}, (X \sm S')_{k_S})$. By assumption, $X\sm S$ is a $K(\pi,1)$, hence   $H^i_{\et}((X \sm S)_{k_S})=0$ for $i \geq 1$ by condition (iii) of Proposition~\ref{kpi1cond}. This implies isomorphisms
\[
H^i_{\et}\big((X \sm S')_{k_S}\big) \stackrel{\sim}{\to} \ressum_{{\mathfrak p} \in S'\sm S (k_S)} H^{i+1}_{\mathfrak p}\big(((X \sm S)_{k_S})_{\mathfrak p}\big)
\]
for $i\geq 1$. As $k_S$ realizes the maximal unramified $p$-extension of $k_\q$ for all $\q\in S'\sm S$, the schemes $((X \sm S)_{k_S})_{\mathfrak p}$, $\p\in S'\sm S (k_S) $ have trivial cohomology with values in $\Z/p\Z$ and we obtain isomorphisms
\[
H^i((k_S)_\p)\stackrel{\sim}{\to} H^{i+1}_{\mathfrak p}\big(((X \sm S)_{k_S})_{\mathfrak p}\big)
\]
for $i\geq 1$. These groups vanish for $i\geq 2$. This implies
\[
H^i_{\et}((X\sm S')_{k_S})=0
\]
for $i\geq 2$. The scheme $(X\sm S')_{k_{S'}}$ is the universal pro-$p$ covering of
$(X\sm S')_{k_{S}}$.  The Hochschild-Serre spectral sequence yields an inclusion
\[
H^2(\Gal(k_{S'}|k_S)) \hookrightarrow H^2_{\et}((X\sm S')_{k_{S}})=0.
\]
Hence $\Gal(k_{S'}|k_S)$ is a free pro-$p$-group and

\[
H^1(\Gal(k_{S'}|k_S)) \stackrel{\sim}{\to} H^1_{\et}((X\sm S')_{k_S})
\cong \ressum_{{\mathfrak p}\in S'\sm S(k_S)} H^1((k_S)_{\mathfrak p}).
\]
This shows that each ${\mathfrak p} \in S'\sm S(k_S)$ ramifies in $k_{S'}|k_S$, and since the Galois group is free, $k_{S'}$ realizes the maximal $p$-extension of $(k_S)_{\mathfrak p}$. In particular,
\[
H^1(T_{\mathfrak p}(k_{S'}|k_S)) \cong H^1((k_S)_{\mathfrak p})
\]
for all $\mathfrak p \in S'\sm S (k_S)$.
Using that $\Gal(k_{S'}|k_S)$ is free, the Hochschild-Serre spectral sequence induces an isomorphism
\[
0=H^2_\et((X\sm S')_{k_{S}})\mapr{\sim} H^2_\et((X\sm S')_{k_{S'}})^{\Gal(k_{S'}|k_S)}.
\]
Hence $H^2_\et((X\sm S')_{k_{S'}})=0$, since $\Gal(k_{S'}|k_S)$ is a pro-$p$-group.  Condition (iii) of Proposition~\ref{kpi1cond}  implies that $X\sm S'$ is a $K(\pi,1)$.
\end{proof}

\begin{corollary} \label{riemann}
Assume that $X\sm S$ is a $K(\pi,1)$, and let $S\subset S'$ be a finite set of primes of norm $\equiv 1 \bmod p$. Assume that each $\q\in S'\sm S$ does not split completely in $k_S$. Then the arithmetic form of Riemann's existence theorem holds, i.e.\ the natural homomorphism
\[
\freeproductmed_{{\mathfrak p} \in S'\backslash S(k_S)} T_\p(k_{S'}|k_S) \longrightarrow \Gal(k_{S'}|k_S)
\]
is an isomorphism. Here $T_\p$ is the inertia group and $\freeproductmed$ denotes the free pro-$p$-product of a bundle of pro-$p$-groups, cf.\ \cite{NSW}, Ch.\,IV,\,\S3.
\end{corollary}

\begin{proof} By Proposition~\ref{enlarge} and by the calculation of the cohomology of a free product (\cite{NSW}, 4.3.10 and 4.1.4), $\phi$ is a homomorphism between free pro-$p$-groups which induces an isomorphism on mod~$p$ cohomology. Therefore $\phi$ is an isomorphism.
\end{proof}

\section{Proof of Theorems 3 and 5}
\begin{theorem} \label{allesda} Assume that $X\sm S$ is a $K(\pi,1)$ and $G_S\neq 1$. Then $k_S$ realizes the maximal  $p$-extension $k_\p (p)$ of the local field $k_\p$ for all $\p\in S$.
\end{theorem}

\begin{proof} The decomposition groups of primes in $S$ have infinite index by Corollary~\ref{infinite_index}. By Corollary~\ref{finite_ext}, we may replace $k$ by a finite subextension in $k_S$, and therefore assume that $\#S\geq 2$.

By Proposition~\ref{wieweitlokal}, it suffices to show that each $\p \in S$ ramifies in $k_S$. Let $\p \in S$ be a prime which does not ramify in $k_S$ and put $S'=S\sm \{\p\}$. By Lemma~\ref{whenramifies}, the natural injection $\phi\colon\ V_{S} \hookrightarrow V_{S'}$ is not an isomorphism. By Proposition~\ref{VSchange}, the cokernel of $\phi$ is one-dimensional. By Theorem~\ref{globcoh}, we obtain
\[
h^2(X \sm S')= h^2(X \sm S).
\]
As $G_S=G_{S'}$, we have $\cd G_{S'}=2$ and
\[
h^2(G_{S})= h^2(G_{S'}) \leq h^2(X \sm S') = h^2(X \sm S).
\]
As $X\sm S$ is a $K(\pi,1)$, equality holds. Therefore the injection $H^2(G_{S'})\hookrightarrow H^2_\et(X\sm S')$ is an isomorphism. By Corollary~\ref{condforkpi1}, $X\sm S'$ is a $K(\pi,1)$. By Proposition~\ref{wieweitlokal}, $\p$ does not split completely in $k_{S'}=k_S$. By Proposition~\ref{enlarge}, $k_S$ realizes the maximal $p$-extension of $k_\p$. This yields a contradiction.
\end{proof}

Now we are in the position to show Theorem~\ref{1.4}.

\begin{proof}[Proof of Theorem~\ref{1.4}]
We have $H^2_\et((X\sm S)_{k_S})=0$ by condition (iii) of Proposition~\ref{kpi1cond}. By Theorem~\ref{allesda}, the local cohomology groups $H^2_\p((X_{k_S})_\p)$ vanish for all $\p \in S(k_S)$. Therefore the excision sequence yields $H^2_\et(X_{k_S})=0$. By the flat duality theorem of Artin-Mazur (\cite{Mi}, III Corollary 3.2) we have
$H^2_\et(X_K)^\vee\cong H^1_\fl (X_K,\mu_p)$ for each finite subextension $K$ of $k$ in $k_S$. Hence
\[
\varprojlim_{K\subset k_S} H^1_\fl(X_K,\mu_p)=0.
\]
The flat Kummer sequence $0\to \mu_p \to \G_m \stackrel{\cdot p}{\to} \G_m \to 0$ implies compatible exact sequences
\[
0 \to \O_K^\times /p \to H^1_\fl(X_K,\mu_p) \to \null_p H^1_\fl(X_K,\G_m)
\]
for all $K$.  We obtain
\[
\varprojlim_{K\subset k_S} \O_K^\times /p =0\, .
\]
The topological Nakayama-Lemma (see \cite{NSW}, Corollary 5.2.8) for the compact $\Z_p$-module $\varprojlim \O_K^\times \otimes \Z_p$ therefore implies
\[
\varprojlim_{K\subset k_S} \O_K^\times \otimes \Z_p =0\, .
\]
Tensoring the exact sequences (cf.\ \cite{NSW}, Lemma 10.3.11)
\[
0 \to \O_K^\times \to \O_{K,S}^\times \to \bigoplus_{\p\in S(K)} (K_\p^\times/U_\p) \to \Cl(K) \to \Cl_S(K) \to 0
\]
by (the flat $\Z$-algebra) $\Z_p$, we obtain exact sequences of finitely generated, hence compact, $\Z_p$-modules. The field $k_S$ admits no unramified $p$-extensions. Therefore class field theory implies $\varprojlim_K \Cl(K)(p)=0$, where $K$ runs through all finite subextensions of $k$ in $k_S$. Thus, passing to the projective limit over $K$, we obtain the exact sequence
\[
0 \to \varprojlim_{K\subset k_S} \O_K^\times \otimes \Z_p \to \varprojlim_{K\subset k_S} \O_{K,S}^\times \otimes \Z_p \to \varprojlim_{K\subset k_S} \bigoplus_{\p\in S(K)} (K_\p^\times /U_\p) \otimes \Z_p \to 0.
\]
As $k_S$ realizes the maximal unramified $p$-extension of $k_\p$ for all $\p \in S$, local class field theory implies the vanishing of the right hand limit.
Therefore the result for the $S$-units follows from the corresponding result for the units.
\end{proof}

We have proven all assertions but the statement on the dualizing module in Theorem~\ref{1.1}. In \cite{S}, Th.~5.2 we showed this statement under the assumption that $k_S$ realizes the maximal $p$-extension $k_\p(p)$ of $k_\p$ for all $\p\in S$. This assumption has been shown above, hence the result follows.

\vskip2cm

\noindent \footnotesize{Alexander Schmidt, NWF I - Mathematik, Universit\"{a}t Regensburg, D-93040
Regensburg, Deutschland. email: alexander.schmidt@mathematik.uni-regensburg.de}


\begin{thebibliography}{NSW}

\bibitem[AGV]{AGV} M. Artin, A. Grothendieck and J.-L. Verdier {\it
Th\'{e}orie des Topos et Cohomologie \'{E}tale des Sch\'{e}mas}. Lecture Notes in Math. 269, 270, 305,
Springer, Heidelberg, 1972/73.

\bibitem[AM]{AM} M. Artin and B. Mazur {\it \'{E}tale homotopy}. Lecture Notes in Math. No. 100 Springer-Verlag, Berlin-New York 1969

\bibitem[FM]{FM} J.-M. Fontaine, B. Mazur {\it Geometric Galois representations}. In: Elliptic Curves,
 Modular Forms, \& Fermat's Last Theorem, edited by J.\ Coates, and S.T.\ Yau. International Press, Boston 1995

\bibitem[Fr]{Fr} E.~M.~Friedlander {\it Etale homotopy of simplicial schemes}. Princeton University Press 1982.

\bibitem[La]{La} J.~P.~Labute {\it Mild pro-$p$-groups and Galois groups of $p$-extensions of\/ $\Q$}.  J.~Reine und angew.~Math. {\bf 596} (2006), 155--182.

\bibitem[Lz]{Lz} M.~Lazard {\it Sur les groupes nilpotents et les anneaux de Lie.}
Ann.\ Ec.\ Norm.\ Sup.\ {\bf 71} (1954) 101-190

\bibitem[Ne]{Ne} J. Neukirch {\it Algebraic Number Theory}. Grundlehren der math.\ Wissenschaften Bd.\ 322, Springer 1999.

\bibitem[NSW]{NSW}
J.~Neukirch, A.~Schmidt, K.~Wingberg, \textit{Cohomology of Number Fields},
Grund\-lehren der math.\ Wiss.\ Bd.\ 323, Springer-Verlag 2000.

\bibitem[Ma]{Ma} B. Mazur {\it Notes on \'{e}tale cohomology of number fields}.  Ann. Sci. \'{E}cole Norm. Sup. (4)  {\bf 6}  (1973), 521--552.

\bibitem[Mi]{Mi} J.S. Milne {\it Arithmetic duality theorems}. Academic Press 1986.

\bibitem[S1]{S} A.~Schmidt \textit{Circular sets of prime numbers and $p$-extension of the rationals}. J.~Reine und angew.~Math.\ {\bf 596} (2006), 115--130.

\bibitem[S2]{twoinf} A. Schmidt {\it On the relation between $2$ and $\infty$ in Galois cohomology of number fields}. Compositio Math. {\bf 133}  (2002),  no. 3, 267--288.

\bibitem[Vo]{Vo}
D.~Vogel, \textit{Circular sets of primes of imaginary quadratic number fields},
Preprints der Forschergruppe {\it Algebraische Zykel und $L$-Funktionen} Regensburg/Leipzig Nr. 5, 2006.\\
{\tt http://www.mathematik.uni-regensburg.de/FGAlgZyk}

\bibitem[Zi]{Zi} T. Zink {\em Etale cohomology and duality in number fields}. Appendix 2
to K. Haberland: Galois cohomology of algebraic number fields. Dt.\ Verlag der
Wissenschaften, Berlin 1978
\end{thebibliography}
\end{document}